\documentclass[12pt]{amsart}

\usepackage{cite} 
\usepackage{color} 
\usepackage{graphicx} 
\usepackage{latexsym,amssymb,amsmath,amsfonts} 
\usepackage{mathrsfs} 
\usepackage{multirow}
\usepackage{bigstrut}
\usepackage{algorithm}
\usepackage{algorithmic}

\allowdisplaybreaks[4]
\newtheorem{mylemma}{Lemma}[section]
\newtheorem{mydefinition}{Definition}[section]
\newtheorem{mytheorem}{Theorem}[section]
\newtheorem{myproposition}{Proposition}[section]
\newtheorem{mycorollary}{Corollary}[section]
\newtheorem{myexample}{Example}[section]

\definecolor{Light}{gray}{0.85}

\def\abs#1{\left\vert #1 \right\vert}
\def\allpoly{\mbox{$\re\langle X \rangle$}}
\def\allpolyell{\mbox{$\re^{\ell}\langle X \rangle$}}
\def\allpolyx0degn{\mbox{$P_n$}}

\def\allseries{\mbox{$\re\langle\langle X \rangle\rangle$}}

\def\allseriesell{\mbox{$\re^{\ell} \langle\langle X \rangle\rangle$}}

\def\allseriesLC{\mbox{$\re_{LC}\langle\langle X \rangle\rangle$}}

\def\allseriesellLC{\mbox{$\re^{\ell}_{LC}\langle\langle X \rangle\rangle$}}

\def\bull{\rule{0.08in}{0.08in}} 

\def\card{{\rm card}}
\def\charseries{{\rm char}}

\newcommand{\comment}[1]{} 

\def\diag{{\rm diag}}

\def\eqref#1{(\ref{#1})} 

\def\expup{{\rm e}}

\def\lieseries{\widehat{\mathcal L}(X)}

\def\Malcevgroup{\widehat{\mathcal G}(X)}

\def\mbf#1{\hbox{\mathversion{bold}$#1$}} 

\def\norm#1{\left\Vert#1\right\Vert}
\def\notin{{\not\in}}

\def\openbull{\framebox[0.08in][c]{$\;$}} 

\def\re{{\mathbb R}} 

\def\shuffle{{\scriptscriptstyle \;\sqcup \hspace*{-0.05cm}\sqcup\;}}

\def\supp{{\rm supp}}

\def\begals{\[\begin{aligned}}
\def\endals{\end{aligned}\]}
\def\begce{\begin{center}}
\def\endce{\end{center}}
\def\begar{\begin{array}}
\def\endar{\end{array}}
\def\begeq{\begin{equation}}
\def\endeq{\end{equation}}
\def\begdi{\begin{displaymath}}
\def\enddi{\end{displaymath}}
\def\begdis{\begin{eqnarray*}}
\def\enddis{\end{eqnarray*}}
\def\begeqa{\begin{eqnarray}}
\def\endeqa{\end{eqnarray}}
\def\begdes{\begin{description}}
\def\enddes{\end{description}}
\def\begit{\begin{itemize}}
\def\endit{\end{itemize}}
\def\begen{\begin{enumerate}}
\def\enden{\end{enumerate}}
\def\beglar{\left[\begin{array}}
\def\endrar{\end{array}\right]}
\def\begle{\begin{mylemma}}
\def\endle{\end{mylemma}}
\def\begde{\begin{mydefinition}}
\def\endde{\end{mydefinition}}
\def\begth{\begin{mytheorem}}
\def\endth{\end{mytheorem}}
\def\begco{\begin{mycorollary}}
\def\endco{\end{mycorollary}}
\def\begprop{\begin{myproposition}}
\def\endprop{\end{myproposition}}
\def\begex{\begin{myexample} \rm}
\def\endex{\hfill\openbull \end{myexample} \vspace*{0.15in}}
\def\begexer{\begin{myexercise}}
\def\endexer{\end{myexercise}}

\def\begres{\noindent{\bf Remarks}:\begin{enumerate}}
\def\endres{\end{enumerate} \par}
\def\begpr{\noindent{\em Proof:}$\;\;$}
\def\endpr{\hfill\bull \vspace*{0.15in}}
\def\begtab{\begin{tabular}}
\def\endtab{\end{tabular}}
\def\rref#1{(\ref{#1})}
\def\begfi{\begin{myfigure}}
\def\endfi{\end{myfigure}}

\def\shuff#1#2{\mathbin{
      \hbox{\vbox{\hbox{\vrule \hskip#2 \vrule height#1 width 0pt}\hrule}\vbox{\hbox{\vrule \hskip#2 \vrule height#1 width 0pt\vrule }\hrule}}}}
\def\shuffl{{\mathchoice{\shuff{5pt}{3.5pt}}{\shuff{5pt}{3.5pt}}{\shuff{3pt}{2.6pt}}{\shuff{3pt}{2.6pt}}}}
\def\shuffle{{\, \shuffl \,}}

\def\lbracedef{\left\{\begin{array}{@{\hspace*{2pt}}c@{\hspace*{3pt}}c@{\hspace*{3pt}}l}} 
\def\allseriesL{\mbox{$\re[[L]]$}}

\def\Lbasis{\ell} 
\def\snorm#1{{|||#1|||}}
\def\allseriesY{\mbox{$\re\langle\langle Y\rangle\rangle$}}
\def\allseriesXxY{\mbox{$\re\langle\langle X\otimes Y\rangle\rangle$}}
\def\Lpoly1{\re[L\cup\emptyset]}
\def\Lseries1{\re[[L\cup\emptyset]]}
\def\allpolyprop{\mbox{$\re_p\langle X \rangle$}}
\def\allseriesprop{\mbox{$\re_p\langle\langle X \rangle\rangle$}}
\def\compeff{CE}
\def\hatcompeff{\widehat{CE}}

\newcounter{mycount}

\begin{document}

\title [Lyndon Word Transduction for Chen--Fliess Series]
{Lyndon Word Transduction to Improve the \\ Efficiency of Computing Chen--Fliess Series}

\author{W. Steven Gray}
\address{(WSG) Department of Electrical and Computer Engineering, Old Dominion University, Norfolk, Virginia 23529, USA}
\email{sgray@odu.edu}
\urladdr{http://www.ece.odu.edu/~sgray}

\author{Fatima Raza}
\address{(FR) Department of Electrical and Computer Engineering, Old Dominion University, Norfolk, Virginia 23529, USA}
\email{fraza001@odu.edu}

\author{Lance Berlin}
\address{(LB) Department of Electrical and Computer Engineering, Old Dominion University, Norfolk, Virginia 23529, USA}
\email{lberl001@odu.edu}

\author{Luis~A.~Duffaut Espinosa}
\address{(LADE) Department of Electrical and Biomedical Engineering, University of Vermont, Burlington, VT 05405, USA}
\email{lduffaut@uvm.edu}
\urladdr{https://www.uvm.edu/~lduffaut}

\date{}

\begin{abstract}
The class of input-output systems representable as Chen--Fliess series arises often in control theory.
One well known drawback of this representation, however, is that the iterated integrals which appear in these
series are algebraically related by the shuffle product.
This becomes relevant when one wants to numerically evaluate these series in applications,
as this redundancy leads to unnecessary computational expense. The general goal of this paper is
to present a computationally efficient way to evaluate these series by introducing what amounts to a change of
basis for the computation.
The key idea
is to use the fact that
the shuffle algebra on (proper) polynomials over a finite alphabet is isomorphic to the polynomial algebra generated by the Lyndon words over this alphabet.
The iterated integrals indexed by Lyndon words contain all the input information needed to compute the output.
The change of basis is accomplished by applying a transduction, that is,
a linear map between formal power series in different alphabets,
to re-index the Chen--Fliess series in terms of Lyndon monomials.
The method is illustrated using a simulation of a continuously stirred-tank reactor system under a cyber-physical attack.
\end{abstract}

\maketitle

\thispagestyle{empty}

\noindent
\tableofcontents


\vspace{0.1in}


\section{Introduction}

The class of input-output systems representable as Chen--Fliess series arises often in control theory \cite{Fliess_81,Fliess_83,Isidori_95}.
Their convergence properties are well understood \cite{Beauchard-etal_23,Fliess_81,Gray-Wang_02}, and their structure as series of iterated integrals
makes them amenable to many forms of algebraic analysis \cite{Ferfera_79,Ferfera_80,Fliess_81,Fliess_83,Gray-etal_Automatic_14,Gray-etal_SCL14}.
One well known drawback of this representation, however, is that iterated integrals are algebraically related, in the simplest instance by the
integration by parts formula, and more generally in terms of a shuffle algebra over the rational field \cite{Fliess_81,Reutenauer_93}.
Therefore, such series representations contain redundant
information.
This fact becomes relevant when one wants to numerically evaluate these series in applications,
as this redundancy leads to unnecessary computational expense. Therefore, the general goal of this paper is
to present a computationally efficient way to evaluate these series by introducing what amounts to a change of
basis for the computation.

The closest related literature to the proposed problem is that for computing the signature of a path in $\re^m$ \cite{Lyons-etal_07}.
Signatures have found wide application in finance \cite{Bayer-etal_25,Kuzmanovic_18} and machine learning \cite{Chevyrev-Kormilitzin_25,Lyons-McLeod_25}.
Similar to Chen--Fliess series, they involve (formal) sums of iterated integrals, but there are some key differences:
\setcounter{mycount}{0}
\begin{list}{\roman{mycount})}{\usecounter{mycount} \itemsep 0.0in \parsep 0.0in \itemindent -0.1in}
\item Signatures operate on paths, that is, parameterized curves that are normalized to have unit length and unit speed along the path. Chen--Fliess series
are defined on control signals, which are usually not normalized in any way.
\item Signatures (and their close relative {\em Chen series} \cite{Sussmann_86}) are exponential Lie series. Chen--Fliess series are generally not of this type.
\item Chen--Fliess series involve a generating series, namely a formal power series that can be computed from an underlying
state space model (if one exists \cite{Fliess_81,Fliess_83,Isidori_95}) or identified from input-output data \cite{Gray-etal_Auto20,Gray-etal_ACC22}. It can also be completely arbitrary.
The coefficients of the generating series provide weights for the iterated integrals.
Signatures have no such counterpart.
\end{list}
There are a number of software packages available for efficiently computing
signatures \cite{Genet-Inzirillo_25,Kidger-Lyons_21,Lyons-Ni_17,Reizenstein_17,Reizenstein-Graham_20,Tong_24}.
The key idea behind optimizing this computation is to use the fact that the logarithm of a signature is a Lie series in the free Lie algebra
over an alphabet $X$ containing $m$ letters.
There are known basis sets for this Lie algebra, for example, the Lyndon word basis \cite{Lothaire_83}.
In which case, the only iterated integrals that need to be computed are those corresponding to the Lyndon words $L$ over $X$, which are
considerably fewer than in the general case.
In many machine learning
applications, computing the log-signature is in fact the primary goal.
In the context of Chen--Fliess series, however, the log signature is not so
relevant an object. Nevertheless, Lyndon words are still useful in solving the proposed problem.

The main idea behind a more efficient computation of the Chen--Fliess series is to use the fact that
the (commutative) shuffle algebra on the linear space of (proper) polynomials in $X$ is isomorphic to the commutative polynomial algebra
generated by the Lyndon words $L$ over $X$
\cite{Grinberg-Reiner_00,Hazewinkel_01,Lothaire_83,Radford_79,Reutenauer_93}.
The shuffle algebra is the graded dual of the universal enveloping algebra of the free Lie algebra on $X$ \cite{Lothaire_83}.
It turns out that the iterated integrals indexed by Lyndon words contain all the input information needed to compute the output.
The objective therefore is to rewrite the Chen--Fliess series in terms of
an expansion indexed solely by Lyndon word monomials.
This is accomplished by applying a {\em transduction} from the finite alphabet $X$ to the countably infinite alphabet $L$.
A transduction is a linear map between formal power series in different alphabets \cite{Fliess_74a,Jacob_75,Jacob_75diss,Kuich-Salomaa_86,Salomaa-Soittola_78}.
In the current context, it is useful to first present a transduction as an action on the universal control systems of Kawski and Sussmann \cite{Kawski-Sussmann_97}.
These behave in much the same way that linear transformations act on finite dimensional state space realizations, albeit universal control systems are {\em not}
finite dimensional.
This in turn induces a transduction on the Chen--Fliess representation of the input-output map. These {\em Chen--Fliess transductions}
will initially be viewed as formal maps without regard to convergence. But for the particular transduction of interest from $X$ to $L$,
a matrix representation will be derived in order to perform convergence analysis. In particular, it will be shown that this
transduction preserves the local convergence of the original series. It should also be noted that transductions were first used in the context of Chen--Fliess
series by Ferfera in \cite{Ferfera_79} to provide a sufficient condition under which the series interconnection of two rational (i.e., bilinear) systems
remains rational. But the notion of a Chen--Fliess transduction as presented here is only implicit in this earlier work.

Finally, the proposed solution is evaluated using
the concept of {\em computational efficiency}. This is a relative measure of the reduction in the number of iterated integrals
needed to represent a generic input-output map up to a given word length (level in the language of signatures) in the new series expansion in terms of $L$ as compared to the original series in terms of $X$. An asymptotic approximation of this measure as a function of word length is given
in order to easily assess a priori the potential benefits of applying a Chen--Fliess transduction in a given application. A number of examples are considered, including
one involving a type of cyber-physical attack on a continuously stirred-tank reactor (CSTR) system \cite{Gray-etal_CISS21,Cao-etal_20,Jeon-Eun_19,Park-etal_18}.

The paper is organized as follows. In the next section, a summary of the mathematical preliminaries is given to establish the notation and the
key concepts needed for the rest of the paper.
Then, in Section~\ref{sec:Chen-Fliess-transductions}, the notion of a Chen--Fliess transduction is presented.
This is applied in Section~\ref{sec:Lyndon-coordinates} to rewrite Chen--Fliess series in terms of Lyndon words. Applications to improving the computational
efficiency of such series representations are described in Section~\ref{sec:computational-efficiency}. The conclusions of the paper are given in the final section.

\section{Preliminaries}

An {\em alphabet} $X=\{ x_0,x_1,$ $\ldots,x_m\}$ is any nonempty finite set
of symbols referred to as {\em
letters}. A {\em word} $\eta=x_{i_1}\cdots x_{i_k}$ is a finite sequence of letters from $X$.
The number of letters in a word $\eta$, written as $\abs{\eta}$, is called its {\em length}.
The empty word, $\emptyset$, has length zero.
The collection of all words having length $k$ is denoted by
$X^k$. Define $X^{\leq j}=\bigcup_{k=0}^j X^k$, $X^+=\bigcup_{k> 0} X^k$, and $X^\ast=\bigcup_{k\geq 0} X^k$.
The set $X^\ast$ is a monoid under the concatenation product with the empty word acting as the unit.
Any mapping $c:X^\ast\rightarrow
\re^\ell$ is called a {\em formal power series}.
Often $c$ is
written as the formal sum $c=\sum_{\eta\in X^\ast}( c,\eta)\eta$,
where the {\em coefficient} $( c,\eta)\in\re^\ell$ is the image of
$\eta\in X^\ast$ under $c$.
The {\em support} of $c$, $\supp(c)$, is the set of all words having nonzero coefficients.
A series $c$ is called {\em proper} if $\emptyset\notin\supp(c)$.
The set of all noncommutative formal power series over the alphabet $X$ is
denoted by $\allseriesell$. The subset of series with finite support, i.e., polynomials,
is represented by $\allpolyell$.
Each set is an associative $\re$-algebra under the concatenation product and an associative and commutative $\re$-algebra under
the {\em shuffle product}, that is, the bilinear product uniquely specified by the shuffle product of two words
\begdi
	(x_i\eta)\shuffle(x_j\xi)=x_i(\eta\shuffle(x_j\xi))+x_j((x_i\eta)\shuffle \xi),
\enddi
where $x_i,x_j\in X$, $\eta,\xi\in X^\ast$ and with $\eta\shuffle\emptyset=\emptyset\shuffle\eta=\eta$ \cite{Fliess_81}.

\subsection{Chen--Fliess series and universal control systems}

Given any $c\in\allseriesell$, one can associate a causal
$m$-input, $\ell$-output operator, $F_c$, in the following manner.
Let $\mathfrak{p}\ge 1$ and $t_0 < t_1$ be given. For a Lebesgue measurable
function $u: [t_0,t_1] \rightarrow\re^m$, define
$\norm{u}_{\mathfrak{p}}=\max\{\norm{u_i}_{\mathfrak{p}}: \ 1\le
i\le m\}$, where $\norm{u_i}_{\mathfrak{p}}$ is the usual
$L_{\mathfrak{p}}$-norm for a measurable real-valued function,
$u_i$, defined on $[t_0,t_1]$.  Let $L^m_{\mathfrak{p}}[t_0,t_1]$
denote the set of all measurable functions defined on $[t_0,t_1]$
having a finite $\norm{\cdot}_{\mathfrak{p}}$ norm and
$B_{\mathfrak{p}}^m(R)[t_0,t_1]:=\{u\in
L_{\mathfrak{p}}^m[t_0,t_1]:\norm{u}_{\mathfrak{p}}\leq R\}$.
Assume $C[t_0,t_1]$
is the subset of continuous functions in $L_{1}^m[t_0,t_1]$. Define
inductively for each word $\eta=x_i\bar{\eta}\in X^{\ast}$ the map $E_\eta:
L_1^m[t_0, t_1]\rightarrow C[t_0, t_1]$ by setting
$E_\emptyset[u]=1$ and letting
\[E_{x_i\bar{\eta}}[u](t,t_0) =
\int_{t_0}^tu_{i}(\tau)E_{\bar{\eta}}[u](\tau,t_0)\,d\tau, \] where
$x_i\in X$, $\bar{\eta}\in X^{\ast}$, and $u_0=1$. The
{\em Chen--Fliess series} corresponding to $c\in\allseriesell$ is
\begdi
y(t)=F_c[u](t) =
\sum_{\eta\in X^{\ast}} ( c,\eta) \,E_\eta[u](t,t_0) 
\enddi
\hspace*{-0.05in}\cite{Fliess_81,Fliess_83,Isidori_95}.
If there exist real numbers $K,M>0$ such that
\begeq
\abs{( c,\eta)}\le K M^{|\eta|}|\eta|!,\;\; \forall\eta\in X^{\ast},
\label{eq:local-convergence-growth-bound}
\endeq
then $F_c$ constitutes a well defined mapping from
$B_{\mathfrak p}^m(R)[t_0,$ $t_0+\Delta]$ into $B_{\mathfrak
q}^{\ell}(S)[t_0, \, t_0+\Delta]$ for sufficiently small $R,\Delta >0$ and some $S>0$,
where the numbers $\mathfrak{p},\mathfrak{q}\in[1,\infty]$ are
conjugate exponents, i.e., $1/\mathfrak{p}+1/\mathfrak{q}=1$  \cite{Gray-Wang_02}.
(Here, $\abs{z}:=\max_i \abs{z_i}$ when $z\in\re^\ell$.)
Any series $c$ satisfying \rref{eq:local-convergence-growth-bound} is called
{\em locally convergent}.
The set of all
locally convergent series is denoted by $\allseriesellLC$, and $F_c$ is
referred to as a {\em Fliess operator}.

Given any Fliess operator $y=F_c[u]$, one can associate the {\em universal control system}
\begin{subequations} \label{eq:UCS}
\begin{align} 
\dot{z}&=\sum_{i=0}^m x_iu_i z,\;\; z(0)={\mbf 1} \label{eq:UCS-state-equation} \\
y_k&=( c_k,z),\;\; k=1,2,\ldots,\ell, \label{eq:UCS-output-equation}
\end{align}
\end{subequations}
where $z=\sum_{\eta\in X^\ast} \eta\,E_\eta[u]$, $\mbf{1}:=1\emptyset$, and $c_k$ is the $k$-th component series of $c$.
Here $z$ is an exponential Lie series known as the {\em Chen series}. It can be viewed as an element of the formal Lie group of all such series, $\Malcevgroup$ \cite{Gray-etal_SCL21,Kawski-Sussmann_97}.
Analogous to standard Lie group theory, the formal tangent space at the unit $\mathbf{1}$, $T_{\mathbf{1}}\Malcevgroup$,
is identified with the Lie algebra of $\Malcevgroup$, denoted here by $\lieseries$.
For any $p\in\lieseries$, there is a corresponding tangent vector at $z\in\Malcevgroup$
written as the linear functional
\begdi
V_p(z):\allseriesLC\rightarrow\re,\, e\mapsto V_p(z)(e):= (e,pz)
\enddi
and satisfying the Leibniz rule. In \rref{eq:UCS-state-equation}, this corresponds to setting $p=\sum_{i=0}^m x_iu_i$.
In what follows, \rref{eq:UCS} will be abbreviated as $(V_p,\mbf{1},c)$.

One can also view \rref{eq:UCS} as an infinite dimensional state space realization with states $z_\eta:=(z,\eta)$ indexed by words $\eta\in X^\ast$. In which case,
$\dot{z}_{x_i\eta}=z_\eta u_i$, and $y$ is a linear combination of these states. A
linear coordinate transformation on such a realization can be described in terms of
a {\em transduction} as defined next.

\subsection{Transductions}

Let $X$ and $Y$ be two arbitrary alphabets.
Let $\allseries$ and $\allseriesY$ be $\re$-linear spaces of formal power series over $X$ and $Y$, respectively.
Any $\re$-linear mapping
$\tau_{XY}:\allseries\rightarrow\allseriesY$ is called a {\em
transduction}\cite{Fliess_74a,Jacob_75,Jacob_75diss,Kuich-Salomaa_86,Salomaa-Soittola_78}.
It is completely specified by
\begin{equation*}\label{eq:trans1}
\tau_{XY}(\eta)=\sum\limits_{\xi\in
Y^*}(\tau_{XY}(\eta),\xi)\xi,\;\; \forall \eta\in X^*,
\end{equation*}
if for each $c\in\allseries$, the summation
\begdi
(\tau_{XY}(c),\xi)=\sum_{\eta\in X^\ast}(c,\eta)(\tau_{XY}(\eta),\xi)
\enddi
converges for all $\xi\in Y^\ast$.
Otherwise, it is only partially defined on $\allseries$.
One can canonically associate with any transduction $\tau_{XY}$ a series
$\tau\in\allseriesXxY$, namely
\begdi \label{eq:trans2} \nonumber
\tau = \sum\limits_{\eta\in X^*}\eta\otimes\tau_{XY}(\eta)
= \sum\limits_{\eta\in X^\ast,\;\xi\in Y^\ast}(\tau_{XY}(\eta),\xi)\,\eta\otimes\xi.
\enddi
From $\tau$ define the {\em adjoint} transduction $\tau^\ast_{XY}:\allseriesY\rightarrow\allseries$ via
\begin{equation*}
\tau^\ast_{XY}(\xi)=\sum\limits_{\eta\in X^\ast}(\tau_{XY}(\eta),\xi)\eta,\;\; \forall \xi\in Y^*,
\end{equation*}
so that
\begdi
(\tau_{XY}(\eta),\xi)=(\eta,\tau^{\ast}_{XY}(\xi)),\;\;\forall \eta\in X^\ast,\;\xi\in Y^\ast.
\enddi
It is  easily verified that $\tau_{XY}^{\ast\ast}=\tau_{XY}$.
Finally, if $\tau_{XY}$ corresponds to a canonical isomorphism between the two
$\re$-linear spaces $\allseries$ and $\allseriesY$, then the inverse transduction is written as $\tau^{-1}_{XY}$.
Its dual is abbreviated by
\begdi
\tau^{\dag}_{XY}(\eta):=(\tau^{-1}_{XY})^\ast(\eta)=\sum\limits_{\xi\in Y^\ast}(\tau^{-1}_{XY}(\xi),\eta)\xi,\;\; \forall \eta\in X^\ast.
\enddi

\subsection{Lyndon words and the shuffle algebra}
\label{subsec:necklaces}

Given an alphabet $X$ with finite cardinality, $\card(X)$, a
{\em necklace} of length $n$ is an $n$-gon in an oriented plane with $n$ vertices each labeled by a letter in $X$
\cite[Chapter 7]{Reutenauer_93}.
Two necklaces are said to be equivalent when they can be made identical via rotation. The set of words generated by all rotations defines an equivalence class. A necklace is called {\em primitive} if only the trivial rotation leaves it invariant. The number of primitive necklaces
of length $n$ is
\begdi
P(n)= \frac{1}{n}\,\sum_{d|n} \mu(d)\, (\card(X))^{n/d},
\enddi
where $\mu$ denotes the {\em M\"{o}bius function}, and $d|n$ is the set of all divisors $d$ of $n$.

Suppose the letters in $X$ are ordered as $x_0<x_1<\cdots< x_m$ so as to induce a lexicographical ordering on $X^+$. A word $\eta\in X^+$ is called a {\em Lyndon word} if all factorizations $\eta = \xi \nu$ with $\xi,\nu \in X^+$ have the property that $\eta < \nu\xi$.  Each primitive necklace contains a single Lyndon word in its equivalence class \cite[Corollary~7.7]{Reutenauer_93}.
Therefore, counting primitive necklaces of length $n$ is the same as counting Lyndon words of length $n$.

\begex
Suppose $X=\{x_0,x_1\}$ so that $\card(X)=2$. The first few Lyndon words are $L := \{l_i\}_{i\geq0} = \{x_0, x_1, x_0 x_1, x_0^2 x_1, x_0 x_1^2, x_0^3 x_1, x_0^2 x_1^2, x_0 x_1^3, ...\}$, where here the ordering is by increasing word length and then lexicographically among words of the same length ({\em shortlex} ordering).
\endex

It is known from \cite{Radford_79} (see also \cite[Section 6]{Grinberg-Reiner_00} and \cite[Chapter 5]{Lothaire_83}) that the shuffle algebra on
the subspace of proper polynomials $\allpolyprop\subset\allpoly$
is isomorphic to the free polynomial algebra on $\re[L]$.
Thus, there exists an $\re$-linear map ${\mathscr L}:\allpolyprop\rightarrow \re[L]$ such that
\begeq \label{eq:shuffle-isomorphism}
{\mathscr L}(\eta_1\shuffle\eta_2)={\mathscr L}(\eta_1){\mathscr L}(\eta_2),\;\;\forall \eta_i\in X^\ast
\endeq
with ${\mathscr L}(l_i)=l_i$.
(It will be useful in some cases to extend the definition so that ${\mathscr L}(\emptyset)=\emptyset$ in order to treat the nonproper case.)
Put another way, the shuffle algebra on $\allpolyprop$ is freely generated by the set of all Lyndon words \cite[Theorem 6.1]{Reutenauer_93}.
The Chen-Fox-Lyndon factorization of a word $\eta\in X^+$\
is a unique (lexicographically) non-increasing product of Lyndon words so that
\begdi 
\eta=l_{i_1}l_{i_2}\cdots l_{i_n},\;\;l_{i_1}\geq l_{i_2}\geq\cdots\geq l_{i_n},
\enddi%
\hspace*{-0.05in}
\cite{Chen-etal_58,Hazewinkel_01,Lothaire_83,Reutenauer_93}. This factorization can be used to provide
a recursive algorithm for computing ${\mathscr L}$ \cite{Gray-etal_SCL24}.

\subsection{Lyndon word asymptotics}

Let $X$ be an alphabet with finite
cardinality $\card(X)$.
Denote by $L(n)$ the number of Lyndon words of length $n$ and
$L_+(n):=\sum_{k=1}^n L(k)$ the total number of
Lyndon words {\em up to} length $n$. The following lemma describes the asymptotic behavior of
these two integer sequences.

\begle \label{le:Lyndon-word-asymptotics}
For $n\gg 1$,
\begin{align}
L(n)&\sim\frac{1}{n}(\card(X))^{n} \label{eq:Lyndon-word-asymptotics} \\
L_+(n)&\sim \frac{1}{n}(\card(X))^{n+1}. \label{eq:Lyndon-word-plus-asymptotics}
\end{align}
\endle

\begpr
As described in Subsection~\ref{subsec:necklaces}, it is known that
\begdi
L(n)=P(n)=\frac{1}{n}\,(\card(X))^{n} + \frac{1}{n}\,\sum_{d|n,\,d\neq 1} \mu(d)\, (\card(X))^{n/d}.
\enddi
Thus, as $n$ approaches infinity, the approximation \rref{eq:Lyndon-word-asymptotics} holds.
Substituting \rref{eq:Lyndon-word-asymptotics} into the definition of $L_+(n)$ and using the asymptotic
approximation $\sum_{k=1}^n a^k/k\sim a^{n+1}/n$ gives \rref{eq:Lyndon-word-plus-asymptotics}.

\endpr

\begin{figure}[tb]
\begin{center}
\includegraphics[scale=0.55]{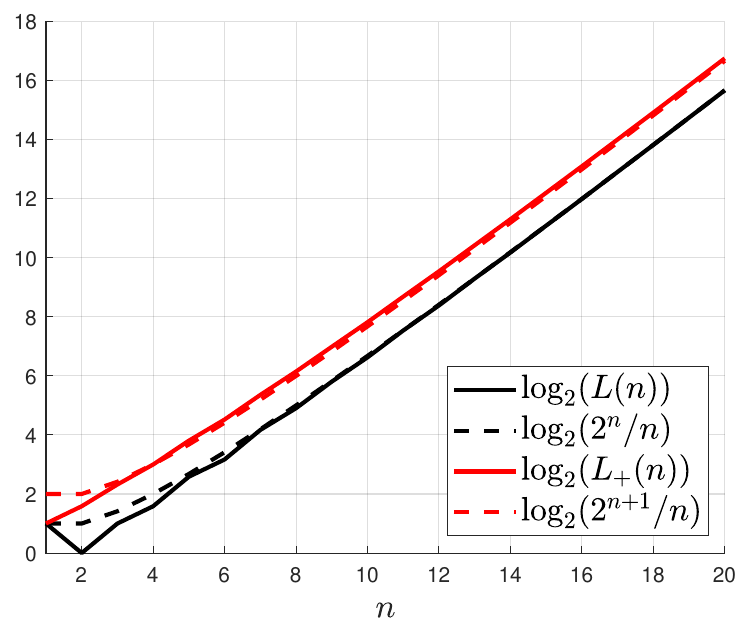}
\caption{Integer sequences $L(n)$ and $L_+(n)$, $n\geq 1$ and their asymptotic approximations when $\card(X)=2$}
\label{fig:FIG1_asymptotics}
\end{center}
\end{figure}

\begex
Suppose $\card(X) = 2$ and $n = 10$. All possible divisors of $10$ are $d = 1,2,5,10$, so that
the total number of Lyndon words of length $10$ is
\begin{align*}
L(10)&=\frac{1}{10}\sum_{d =1,2,5,10} \mu(d) 2^{10/d} \\
&=\dfrac{1}{10}(\mu(1)2^{10} + \mu(2)2^{5} + \mu(5)2^{2} + \mu(10)2) \\
&=\dfrac{1}{10}(2^{10} - 2^{5} - 2^{2} +2) = 99.
\end{align*}
Similarly,
\begdi
\{L(n):n\geq 1\} = \{ 2, 1, 2, 3, 6, 9, 18, 30, 56, 99, ...\},
\enddi
which corresponds to the OEIS integer sequence A001037 \cite{OEIS}.
The integer sequence
\begdi
\{L_+(n):n\geq 1\}=\{2, 3, 5, 8, 14, 23, 41, 71, 127, 226,\ldots\}
\enddi
is OEIS sequence A062692.
Figure~\ref{fig:FIG1_asymptotics} compares these sequences against their asymptotic
approximations in Lemma~\ref{le:Lyndon-word-asymptotics}.
\endex

\section{Transductions for Universal Control Systems and Chen--Fliess Series}
\label{sec:Chen-Fliess-transductions}

In this section, the action of a transduction applied to a
universal control system and a Chen--Fliess series is presented.
The following theorem describes a {\em universal control system transduction}.

\begth \label{th:CF-transduction-UCS}
Assume $\allseries$ and $\allseriesY$ are $\re$-linear spaces over alphabets $X$ and $Y$, respectively,
which are related by a vector space isomorphism $\tau_{XY}$.
Let $(V_p,\mbf{1},c)$ be a universal control system as given in \rref{eq:UCS}. Then
the action of $\tau_{XY}$ on $(V_p,\mbf{1},c)$ is given by
\begdi
(V_p,\mbf{1},c)\stackrel{\tau_{XY}}{\xrightarrow{\hspace*{1cm}}} (V_{\bar{p}},\bar{z}_0,\bar{c}),
\enddi
where
\begin{align*}
\bar{p}(\bar{z})&= \sum_{i=0}^m \tau_{XY}(x_i\tau^{-1}_{XY}(\bar{z}))u_i\\
\bar{z}&=\tau_{XY}(z) \\
\bar{z}_0&=\tau_{XY}(\mbf{1}) \\
\bar{c}&=\tau^\dag_{XY}(c).
\end{align*}
\endth

\begpr
The expressions for $\bar{p}$, $\bar{z}$, and $\bar{z}_0$ follow directly from applying $\tau_{XY}$ to \rref{eq:UCS-state-equation},
the state $z$, and the initial state $z(0)=\mbf{1}$, respectively.
Likewise, the expression for $\bar{c}$ follows from the definition of $\bar{z}$ and \rref{eq:UCS-output-equation}.
\endpr

The next theorem describes a {\em Chen--Fliess transduction}.

\begth \label{th:CF-transduction}
Assume $\allseries$ and $\allseriesY$ are $\re$-linear spaces over alphabets $X$ and $Y$, respectively,
which are related by a vector space isomorphism $\tau_{XY}$. If $c\in\allseries$, then
\begin{align*}
F_c[u]&=\sum_{\eta\in X^\ast} (c,\eta) E_\eta[u] \\
&= F_{\bar{c}}[u]:=\sum_{\xi\in Y^\ast} (\tau^\dag_{XY}(c),\xi)F_{\tau^\ast_{XY}(\xi)}[u].
\end{align*}
\endth

\begpr
Let $(V_p,\mbf{1},c)$ be the universal control system \rref{eq:UCS} corresponding to $F_c$.
The claim is equivalent to saying in the context of Theorem~\ref{th:CF-transduction-UCS} that
the input-output maps for $(V_p,\mbf{1},c)$ and $(V_{\bar{p}},\bar{z}_0,\bar{c})$ coincide on some neighborhood $W$
of $\mbf{1}$ and $\tau_{XY}(W)$ of $\bar{z}_0$, respectively, where $\bar{z}=\tau_{XY}(z)$.
That is, $y=(c,z)=(\bar{c},\bar{z})$ on their respective neighborhoods.
First observe that
\begin{align*}
\bar{z}&:=\tau_{XY}(z) \\
&=\sum_{\eta\in X^\ast} \tau_{XY}(\eta)E_{\eta}[u] \\
&=\sum_{\eta\in X^\ast} \left[\sum_{\xi\in Y^\ast} (\tau_{XY}(\eta),\xi)\xi\right]E_{\eta}[u] \\
&=\sum_{\xi\in Y^\ast} \xi\left[\sum_{\eta\in X^\ast}(\tau_{XY}(\eta),\xi) E_{\eta}[u]\right] \\
&=\sum_{\xi\in Y^\ast} \xi F_{\sum_{\eta\in X^\ast}(\tau_{XY}(\eta),\xi)\eta}[u] \\
&=\sum_{\xi\in Y^\ast} \xi F_{\tau_{XY}^\ast(\xi)}[u].
\end{align*}
Therefore,
\begin{align*}
y&=(c,z) \\
&=(c,\tau^{-1}_{XY}\tau_{XY}(z)) \\
&=(c,\tau^{-1}_{XY}(\bar{z})) \\
&=\sum_{\eta\in X^\ast} (c,\eta)(\tau^{-1}_{XY}(\bar{z}),\eta).
\end{align*}
Next note that
\begin{align*}
\tau^{-1}_{XY}(\bar{z})&=\tau^{-1}_{XY}\left(\sum_{\nu\in Y^\ast} (\bar{z},\nu)\nu\right) \\
&=\sum_{\nu\in Y^\ast} (\bar{z},\nu)\tau^{-1}_{XY}(\nu),
\end{align*}
so that for any $\eta\in X^\ast$
\begdi
(\tau^{-1}_{XY}(\bar{z}),\eta)=\sum_{\nu\in Y^\ast} (\bar{z},\nu)(\tau^{-1}_{XY}(\nu),\eta).
\enddi
Hence,
\begin{align*}
y&=\sum_{\eta\in X^\ast} (c,\eta)\left[\sum_{\nu\in Y^\ast} (\bar{z},\nu)(\tau^{-1}_{XY}(\nu),\eta)\right] \\
&=\sum_{\nu\in Y^\ast}\left[\sum_{\eta\in X^\ast} (c,\eta)(\tau^{-1}_{XY}(\nu),\eta)\right] (\bar{z},\nu) \\
&=\sum_{\nu\in Y^\ast} (\tau^\dag_{XY}(c),\nu)(\bar{z},\nu) \\
&=(\tau^\dag_{XY}(c),\bar{z}) \\
&=(\bar{c},\bar{z})
\end{align*}
as desired.
\endpr

\section{Lyndon Word Transductions for Chen--Fliess Series}
\label{sec:Lyndon-coordinates}

In this section, Chen--Fliess transductions are applied to rewrite
a given Chen--Fliess series in terms of Lyndon words.
That is, $\allseriesY$ in Theorem~\ref{th:CF-transduction} is taken to be the linear space $\re[L]$.
There is clearly more than one way to do this.
Each such map will be called a {\em Lyndon word transduction}.

First, the mappings ${\mathscr L}$ and ${\mathscr L}^{-1}$ are
described in terms of transductions without
regard to convergence issues.
Then the focus turns to two possible choices for Lyndon word transductions.
The Lyndon word transduction of particular interest has the special property that it can be used
to improve the efficiency of computing Chen--Fliess series as described in the next section.
In preparation for the convergence analysis, matrix representations of
${\mathscr L}$ and ${\mathscr L}^{-1}$ are then presented. Finally, the convergence properties of
Chen--Fliess series under the Lyndon word transduction of interest are developed.

\subsection{Mappings ${\mathscr L}$ and ${\mathscr L}^{-1}$ as transductions}

The goal is to prove that ${\mathscr L}:\allseriesprop\rightarrow \re[[L]]$,
where $\allseriesprop\subset\allseries$ is the subspace of all proper series,
is a graded $\re$-vector space isomorphism.
Hence, ${\mathscr L}$ and ${\mathscr L}^{-1}$ can be viewed
as transductions.
In light of the discussion above, this result is a modest refinement of what is known in the literature.
But the constructive approach used here will be useful later for providing explicit matrix
representations of ${\mathscr L}$ and its inverse.
Henceforth, the distinction will be made between a Lyndon word $(l_i)$ as an
element in $X^+$ and the element $l_i\in L$. That is, in general $(l_i)(l_j)\neq (l_j)(l_i)$,
while $l_il_j=l_jl_i$.
The polynomial case is addressed first in the following theorem.

\begth \label{th:L-isomorphism}
The $\re$-linear map ${\mathscr L}:\allpolyprop\rightarrow \re[L]$
is a graded $\re$-vector space isomorphism.
\endth

\begpr
Let $\Psi:\re[L]\rightarrow \allpolyprop$ be the $\re$-linear map uniquely determined by
\begdi
\Psi:L^+\rightarrow \allpoly,\,l=l_{i_1}\cdots l_{i_n}\mapsto (l_{i_1})\shuffle\cdots\shuffle(l_{i_n}).
\enddi
In light of \rref{eq:shuffle-isomorphism}, it follows that
\begin{align*}
({\mathscr L}\circ \Psi)(l)&=({\mathscr L}\circ \Psi)(l_{i_1}\cdots l_{i_n}) \\
&={\mathscr L}( (l_{i_1})\shuffle\cdots\shuffle(l_{i_n})) \\
&=l_{i_1}\cdots l_{i_n} \\
&=l.
\end{align*}
It is well known
that any word $\eta\in X^+$ can be written in terms of a unique linear combination
of shuffle products of Lyndon words \cite{Lothaire_83,Radford_79}.
That is, there exists a unique $q_\eta\in\re[L]$ such that
\begdi
\eta=\sum_{l\in L^+} (q_\eta,l) \Psi(l).
\enddi
Therefore,
\begin{align*}
(\Psi\circ{\mathscr L})(\eta)&= (\Psi\circ{\mathscr L})\left(\sum_{l\in L^+} (q_\eta,l) \Psi(l)\right)\\
   &= \Psi\left(\sum_{l\in L^+} (q_\eta,l)l\right) \\
   &= \sum_{l\in L^+} (q_\eta,l)\Psi(l) \\
   &=\eta.
\end{align*}
That is, ${\mathscr L}$ is injective with ${\mathscr L}^{-1}=\Psi$.
One cannot immediately deduce from this that ${\mathscr L}$ is surjective since neither vector space $\re[L]$ nor $\allpolyprop$ is finite dimensional.
However, these vector spaces are graded in terms of word length with respect to the alphabet $X$. That is, the degree of
$\eta\in X^+$ is defined to be $\abs{\eta}$, and the degree of $l\in L^+$ is defined by
$\abs{l}=\abs{l_{i_1}\cdots l_{i_k}}=\abs{(l_{i_1})}+\cdots +\abs{(l_{i_k})}$.
Let $\allpolyprop=\bigoplus_{k>0} \re \langle X^k\rangle$ be the corresponding grading of $\allpolyprop$,
where $\re \langle X^k\rangle$ denotes the $k$-dimensional $\re$-vector subspace spanned by the
words of degree $k$. Likewise, $\re[L]=\bigoplus_{k> 0} \re[L^k]$ is the corresponding grading of $\re[L]$,
where $\re[L^k]$ is a $k$-dimensional subspace of $\re[L]$.
A property of the shuffle product yields the key fact that
${\mathscr L}$ is degree preserving on the homogeneous components of these graded spaces.
Therefore, ${\mathscr L}$
restricted to any homogeneous component of $\allpolyprop$, say
\begdi
{\mathscr L}_k: \re\langle X^k\rangle\rightarrow \re[L^k],\;\; k>0
\enddi
is injective. Since $\re\langle X^k\rangle$ is finite dimensional, it is also surjective.
Thus, ${\mathscr L}$ is surjective on the entire vector space $\allpolyprop$.
\endpr

The series case is addressed next.

\begth \label{th:L-isomorphism-series}
The $\re$-linear map ${\mathscr L}:\allseriesprop\rightarrow \re[[L]]$
is a graded $\re$-vector space isomorphism.
\endth

\begpr
The result follows from two key facts. First, the gradings on $\allpolyprop$ and $\re[L]$ extend naturally to gradings on
$\allseriesprop$ and $\re[[L]]$, respectively. Second, the images of words under ${\mathscr L}$ and ${\mathscr L}^{-1}$ as defined above are
always polynomials. This implies that the maps ${\mathscr L}$ and ${\mathscr L}^{-1}$ when extended to formal power series yield series
${\mathscr L}(c)$ and ${\mathscr L}^{-1}(c_L)$ that are locally finite, that is, every coefficient of these series is computed from a finite sum.
For example, given any $c\in\allseriesprop$ and $l\in L^k$ and using the fact that ${\mathscr L}$ is degree preserving, it follows that
\begdi
({\mathscr L}(c),l)=({\mathscr L}_k(c),l)=\sum_{\eta\in X^k} ({\mathscr L}_k(\eta),l).
\enddi
Thus, the mappings ${\mathscr L}:\allseriesprop\rightarrow \re[[L]]$ and ${\mathscr L}^{-1}:\allseries[[L]]\rightarrow \allseriesprop$ are at least well defined.
Similarly, ${\mathscr L}$ is bijective with ${\mathscr L}^{-1}$ being its inverse since it follows from the proof of Theorem~\ref{th:L-isomorphism}
that
\begeq \label{eq:block-defined-L-inverse}
{\mathscr L}_k^{-1} {\mathscr L}_k={\mathscr L}_k {\mathscr L}_k^{-1}=I_k,\;\;\forall k> 0,
\endeq
where $I_k$ denotes the identity map on $\re\langle X^k\rangle$ and $\re[L^k]$.
(This fact will become more obvious when matrix representations of ${\mathscr L}$ and ${\mathscr L}^{-1}$ are presented in
Subsection~\ref{subsec:matrix-representations-L-invL}.)
\endpr

\subsection{Lyndon word transductions}

There are two immediate choices for Lyndon word transductions
$\tau_{XL}$, namely, ${\mathscr L}$ and ${\mathscr L}^\dag$. For reasons that will
become clear shortly, the latter is the more useful choice for the application
of improving the efficiency of computing Chen--Fliess series.
It is also helpful at this point to formally add the empty word to both $X^+$ and $L^+$
with ${\mathscr L}(\emptyset)={\mathscr L}^{-1}(\emptyset)=\emptyset$ so that ${\mathscr L}^{-1}(c)=(c,\emptyset)\emptyset+{\mathscr L}^{-1}(c_p)$,
where $c_p:=c-(c,\emptyset)\emptyset\in\allseriesprop$. The set of zero degree polynomials is denoted by $\re\emptyset$.

\begth \label{th:Lyndon-word-transduction}
Let $c\in\allseries$.
Applying the Chen--Fliess transduction $\tau_{XL}={\mathscr L}^\dag:\allseries\rightarrow \re[[L]]\cup\re\emptyset$
gives
\begin{align*}
F_c[u]&=F_{{\mathscr L}(c)}[u] \\
&=\sum_{l\in L^\ast} (\mathscr L(c),l)E_{{\mathscr L}^{-1}(l)}[u],
\end{align*}
where (by a slight abuse of notation)
\begeq \label{eq:ELinv}
E_{{\mathscr L}^{-1}(l)}=E_{{\mathscr L}^{-1}(l_{i_1}\cdots l_{i_k})}:=E_{(l_{i_1})}\cdots E_{(l_{i_k})}.
\endeq
\endth

\begpr
The identity is a direct consequence of Theorem~\ref{th:CF-transduction}.
In particular, observe that
\begin{align*}
\bar{c}&:=\tau_{XL}^\dag(c)=({\mathscr L}^\dag)^\dag(c)={\mathscr L}(c) \\
E_{\tau_{XL}^\ast(l)}[u]&=E_{({\mathscr L}^\dag)^\ast(l)}[u]=E_{{\mathscr L}^{-1}(l)}[u].
\end{align*}
\endpr

A key feature of the transduction ${\mathscr L}^\dag$ is the fact that
the identity \rref{eq:ELinv} will be exploited in Section~\ref{sec:computational-efficiency} to reduce the number
of integrals that need to be computed when evaluating a Chen--Fliess series numerically. The
Chen--Fliess transduction ${\mathscr L}$ does
not have this feature. Consider the following simple example.

\begex
If $c=x_1x_0+x_0x_1$, then
computing $F_c[u]=E_{x_1x_0}[u]+E_{x_0x_1}[u]$ would require computing four integrals.
Observe, however, that
\begin{align*}
E_{x_1x_0}[u]+E_{x_0x_1}[u]&=F_{x_1\shuffle x_0}[u] \\
&=E_{(l_1)}[u]E_{(l_0)}[u] \\
&=E_{{\mathscr L}^{-1}(l_1l_0)}[u],
\end{align*}
which is just a reformulation of the integration by parts formula.
So $F_c[u]$ can be computed from only two integrals, namely, $E_{x_0}[u]$ and $E_{x_1}[u]$.
\endex

\begex
Given any (finite) subset $A\subseteq X^\ast$, define its characteristic (polynomial) series to be $\charseries(A)=\sum_{\eta\in A}\eta$.
Likewise, if $B\subseteq L^\ast$, then $\charseries(B)=\sum_{l\in B} l$. In addition,
define ${\mathscr L}(A)=\{{\mathscr L}(\eta):\eta\in A\}$.
Suppose $X=\{x_0,x_1\}$. Observe
\begin{align}
{\mathscr L}(\charseries(X^2))
&={\mathscr L}(x_0x_0+x_0x_1+x_1x_0+x_1x_1) \nonumber \\
&={\mathscr L}(x_0x_0)+{\mathscr L}(x_0x_1)+{\mathscr L}(x_1x_0)+{\mathscr L}(x_1x_1) \nonumber \\
&=\left(\frac{1}{2}l_0l_0\right)+l_2+(-l_2+l_1l_0)+\left(\frac{1}{2}l_1l_1\right) \nonumber \\
&=(l_0+l_1)^2\frac{1}{2!} \nonumber \\
&=(\charseries({\mathscr L}(X)))^2\frac{1}{2!}. \label{ex:LcharX-identity}
\end{align}
Therefore, if $c=\charseries(X^2)$, then
\begin{align*}
F_{c}[u]&=F_{{\mathscr L}(\charseries(X^2))}[u] \\
&=\frac{1}{k!}F_{(\charseries(L^1))^2}[u] \\
&=\frac{1}{k!}F_{(\charseries(L^1))}^2[u].
\end{align*}
That is, $F_c[u]$ is a polynomial in the integrals $E_{x_0}[u]$, $E_{x_1}[u]$.
(The generalization of \rref{ex:LcharX-identity} will be given shortly in Lemma~\ref{le:Tk-growth}.)
\endex

\begin{table}[tb]
\renewcommand{\arraystretch}{1.5}
\begin{center}
\caption{Lyndon word notation summary}
\label{tbl:notation summary}
\begin{tabular}{|c||c|} \hline
symbol & definition \\ \hline\hline
$l_i$ & Lyndon word in $L$ \\ \hline
$l$ & Lyndon monomial in $L^+$ \\ \hline
$\ell_i$ & order basis element for $\re[L]$ \\ \hline
\end{tabular}
\end{center}
\end{table}

\begin{table}[tb]
\renewcommand{\arraystretch}{1.5}
\begin{center}
\caption{First few Lyndon transformation matrices when $X=\{x_0,x_1\}$}
\label{tbl:Lyndon transformation polynomials}
{\scriptsize
\begin{tabular}{|c||c|c|c|c|} \hline
$k$ & $i$ & $\eta_i\in X^k$ & ${\mathscr L}(\eta_i)$ & $T_k$ matrix\\ \hline\hline
0 & 1 & $\emptyset$ & $\emptyset$ & [1] \\ \hline\hline
\multirow{2}[0]*{1} & 2 & $x_0$ & $l_0$ & \multirow{2}[0]*
{
$\left[ \hspace*{-4pt}
\begin{tabular}{*{2}{c@{\hspace*{3pt}}}}
1 & 0 \\
0 & 1
\end{tabular} \hspace*{-1pt}
\right]$} \bigstrut[tb]
\\ \cline{2-4}
& 3 & $x_1$ & $l_1$ & \\ \hline\hline
\multirow{4}[0]*{2}& 4 &$x_0x_0$ & $\frac{1}{2}l_0l_0$ & \multirow{4}[0]*
{
$\left[ \hspace*{-4pt}
\begin{tabular}{*{4}{r@{\hspace*{3pt}}}}
 $\frac{1}{2}$ & 0 &  0 & 0 \\
         0   & 1 & -1 & 0 \\
         0   & 0 &  1 & 0 \\
         0   & 0 &  0 & $\frac{1}{2}$\\
\end{tabular} \hspace*{-1pt}
\right]$} \bigstrut[tb]
\\ \cline{2-4}
& 5 &$x_0x_1$ & $l_2$ &  \\ \cline{2-4}
& 6 &$x_1x_0$ & $-l_2+l_1l_0$ &  \\ \cline{2-4}
& 7 &$x_1x_1$ & $\frac{1}{2}l_1l_1$ &  \\ \hline\hline
\multirow{8}[0]*{3}& 8 & $x_0x_0x_0$ & $\frac{1}{6}l_0l_0l_0$ & \multirow{8}[0]*
{
$\left[  \hspace*{-4pt}
\begin{tabular}{*{8}{c@{\hspace*{3pt}}}}
$\frac{1}{6}$ &  0 &   0  &  0 &   0             &   0      &  0             &  0 \\
         0  &  1 &  -2  &  0 &   1             &   0      &  0             &  0 \\
         0  &  0 &   1  &  0 &  -1             &   0      &  0             &  0 \\
         0  &  0 &   0  &  1 &   0             &  -2      &  1             &  0 \\
         0  &  0 &   0  &  0 &   $\frac{1}{2}$ &   0      &  0             &  0 \\
         0  &  0 &   0  &  0 &   0             &   1      & -1             &  0 \\
         0  &  0 &   0  &  0 &   0             &   0      &  $\frac{1}{2}$ &  0 \\
         0  &  0 &   0  &  0 &   0             &   0      &  0             &  $\frac{1}{6}$ \\
\end{tabular}  \hspace*{-1pt}
\right]$} \bigstrut[tb]
 \\ \cline{2-4}
& 9 &$x_0x_0x_1$ & $l_3$  & \\ \cline{2-4}
& 10 & $x_0x_1x_0$ & $-2l_3+l_2l_0$ & \\ \cline{2-4}
& 11 & $x_0x_1x_1$ & $l_4$ & \\ \cline{2-4}
& 12 & $x_1x_0x_0$ & $l_3-l_2l_0+\frac{1}{2}l_1l_0l_0$ & \\ \cline{2-4}
& 13 & $x_1x_0x_1$ & $-2l_4+l_1l_2$ & \\ \cline{2-4}
& 14 & $x_1x_1x_0$ & $l_4-l_1l_2+\frac{1}{2}l_1l_1l_0$ &  \\ \cline{2-4}
& 15 & $x_1x_1x_1$ & $\frac{1}{6}l_1l_1l_1$ &  \\ \hline
\end{tabular}
}
\end{center}
\end{table}

\begin{table}[tb]
\renewcommand{\arraystretch}{1.5}
\begin{center}
\caption{First few inverse Lyndon transformation matrices when $X=\{x_0,x_1\}$}
\label{tbl:inverse_Lyndon transformation polynomials}
{\scriptsize
\begin{tabular}{|c|c||c|c|c|} \hline
 $k$ & $i$ & $\Lbasis_i\in L^k$ & ${\mathscr L}^{-1}(\Lbasis_i)$ & $T^{-1}_k$ matrix\\ \hline\hline
0 & 1 & $\emptyset$ & $\emptyset$ & [1] \\ \hline\hline
\multirow{2}[0]*{1} & 2 & $l_0$ & $x_0$ & \multirow{2}[0]*
{
$\left[ \hspace*{-4pt}
\begin{tabular}{*{2}{c@{\hspace*{3pt}}}}
1 & 0 \\
0 & 1
\end{tabular} \hspace*{-1pt}
\right]$} \bigstrut[tb]
\\ \cline{2-4}
& 3 & $l_1$ & $x_1$ & \\ \hline\hline
\multirow{4}[0]*{2}& 4 & $l_0l_0$ & $2x_0x_0$ & \multirow{4}[0]*
{
$\left[ \hspace*{-4pt}
\begin{tabular}{*{4}{r@{\hspace*{3pt}}}}
2 &0 &0 &0 \\
0 &1 &1 &0 \\
0 &0 &1 &0 \\
0 &0 &0 &2
\end{tabular} \hspace*{-1pt}
\right]$} \bigstrut[tb]
\\ \cline{2-4}
& 5 & $l_2$ & $x_0x_1$ &  \\ \cline{2-4}
& 6 & $l_1l_0$ & $x_0x_1+x_1x_0$ &  \\ \cline{2-4}
& 7 & $l_1l_1$ & $2x_1x_1$  &  \\ \hline\hline
\multirow{8}[0]*{3} & 8 & $l_0l_0l_0$ & $6x_0x_0x_0$ & \multirow{8}[0]*
{
$\left[  \hspace*{-4pt}
\begin{tabular}{*{8}{c@{\hspace*{3pt}}}}
6 &0 &0 &0 &0 &0 &0 &0 \\
0 &1 &2 &0 &2 &0 &0 &0 \\
0 &0 &1 &0 &2 &0 &0 &0 \\
0 &0 &0 &1 &0 &2 &2 &0 \\
0 &0 &0 &0 &2 &0 &0 &0 \\
0 &0 &0 &0 &0 &1 &2 &0 \\
0 &0 &0 &0 &0 &0 &2 &0 \\
0 &0 &0 &0 &0 &0 &0 &6
\end{tabular}  \hspace*{-1pt}
\right]$} \bigstrut[tb]
 \\ \cline{2-4}
& 9 & $l_3$ & $x_0x_0x_1$ & \\ \cline{2-4}
& 10 & $l_2l_0$ & $2x_0x_0x_1+x_0x_1x_0$   & \\ \cline{2-4}
& 11 & $l_4$ & $x_0x_1x_1$   & \\ \cline{2-4}
& 12 & $l_1l_0l_0$ & $2x_0x_0x_1+2x_0x_1x_0+2x_1x_0x_0$ & \\ \cline{2-4}
& 13 & $l_1l_2$ & $2x_0x_1x_1+x_1x_0x_1$  & \\ \cline{2-4}
& 14 & $l_1l_1l_0$ & $2x_0x_1x_1+2x_1x_0x_1+2x_1x_1x_0$  &  \\ \cline{2-4}
& 15 & $l_1l_1l_1$ & $6x_1x_1x_1$ &  \\ \hline
\end{tabular}
}
\end{center}
\end{table}

\subsection{Matrix representations for ${\mathscr L}$ and ${\mathscr L}^{-1}$}
\label{subsec:matrix-representations-L-invL}

Denote totally ordered monomial bases for vector spaces $\allseries$ and $\re[L]\cup\re\emptyset$ by
$\{\eta_i\}_{i\geq 1}$ and $\{\Lbasis_i\}_{i\geq 1}$, respectively. (A notation summary is given in Table~\ref{tbl:notation summary}.)
In which case, the matrix representations of ${\mathscr L}$ and ${\mathscr L}^{-1}$, namely, $T$ and $T^{-1}$, respectively,
have components
\begin{align*}
[T]_{\Lbasis_i,\eta_j}&=({\mathscr L}(\eta_j),\Lbasis_i) \\
[T^{-1}]_{\eta_i,\Lbasis_j}&=({\mathscr L}^{-1}(\Lbasis_j),\eta_i).
\end{align*}
In particular, if $\Lbasis_j=l_{j_1}\cdots l_{j_n}$, then
\begdi
[T^{-1}]_{\eta_i,\Lbasis_j}=((l_{j_1})\shuffle\cdots \shuffle (l_{j_n}),\eta_i).
\enddi
In light of the gradings on these vector spaces and \rref{eq:block-defined-L-inverse}, both
${\mathscr L}$ and $\mathscr L^{-1}$ have block diagonal matrix representations
$T=\diag(T_0,T_1,\ldots)$ and $T^{-1}=\diag(T_0^{-1},T_1^{-1},\ldots)$, respectively.
The $T_k$, $k\geq 0$ will be referred to as {\em Lyndon transformation matrices}.

\begex
Set $X=\{x_0,x_1\}$ and select the ordered bases for $\allpoly$ and $\re[L]\cup\re\emptyset$ as shown in
Tables~\ref{tbl:Lyndon transformation polynomials} and \ref{tbl:inverse_Lyndon transformation polynomials},
respectively.
The matrices $T_k$ and $T_k^{-1}$ for $k=0,1,2,3$ are also shown
in these tables.
\endex

Growth bounds for $T_k$ and $T^{-1}_k$ are considered next.
Given a matrix $A\in\re^{n\times n}$, define the
seminorm
\begdi
\snorm{A}_\infty=\max_i \abs{\sum_{j=1}^n A_{ij}}.
\enddi
If the entries $A_{ij}$ are all nonnegative, then the matrix norm
$\norm{A}_\infty=\snorm{A}_\infty$.
Given any sequence of matrices $T_k$, $k\geq 0$ with $T_k\in \re^{n_k\times n_k}$, its
{\em growth bounds} are the triples $(K_+,M_+,s_+)$ and $(K_-,M_-,s_-)$  with $K_\pm,M_\pm>0$ and $s_\pm\in\re$
such that
\begdi
K_-M_-^k(k!)^{s_-}\leq\snorm{T_k}_\infty\leq K_+M_+^k(k!)^{s_+},\;\;\forall k\geq 0.
\enddi

Deriving the growth bounds for the sequence of Lyndon transformation matrices and the sequence of inverses
requires two technical lemmas. Proofs for these lemmas are given
in the appendices.
The first lemma is needed for the growth bounds for  $T_k$, $k\geq 0$.

\begle \label{le:Tk-growth}
Let $X$ be a finite alphabet. For any integer $k\geq 0$,
\begdi
{\mathscr L}(\charseries(X^k))=(\charseries({\mathscr L}(X)))^k\frac{1}{k!}.
\enddi
\endle

\begth \label{th:L-growth}
The matrix representation $T_k$ of the $\re$-linear map ${\mathscr L}_k$ has growth bounds
\begdi
\frac{1}{k!}\leq\snorm{T_k}_{\infty}<\frac{1}{\sqrt{\pi/2}}  2^{k}\frac{1}{k!},\;\; k \geq 1.
\enddi
\endth

\begpr
From Lemma~\ref{le:Tk-growth},
the $\Lbasis_i$-th row sum of $T_k$ is
\begin{align*}
\sum_{\eta_j\in X^k} [T]_{\Lbasis_i,\eta_j}&=\sum_{\eta_j\in X^k}({\mathscr L}(\eta_j),\Lbasis_i),\;\;\ell_i\in L^k \\
&=({\mathscr L}(\charseries(X^k)),\Lbasis_i) \\
&=((\charseries({\mathscr L}(X)))^k,\Lbasis_i)\frac{1}{k!}.
\end{align*}
Therefore, each row sum is a positive integer of the form ${k\choose n}\frac{1}{k!}$, $n\in\{0,1,\ldots,k\}$.
Thus, the maximum row sum  cannot fall below $1/k!$ or exceed ${k\choose \lfloor k/2\rfloor}\frac{1}{k!}$.
The integer sequence $a(k)={k\choose \lfloor k/2\rfloor}$, $k\geq 0$ is OEIS sequence A001405 and is known to grow asymptotically
as $a(k) \sim 2^k / \sqrt{k\pi/2}$. Thus, the upper bound is proved.
\endpr

As a side note, the sequence of matrix norms for the case where $\card(X)=2$ is
\begin{align*} 
\{\|T_k\|_{\infty}, 0\leq k\leq 12\}&=\{1, 1, 2, 4, 8, 36, 104,1140,9608,80956,\\
&\hspace{0.3in}  1334692, 37207684, 7919797428\}
\end{align*}
despite the fact that the entries in these matrices are not all integers.
At present, there is no integer sequence in the OEIS database that matches this data.
In addition, this data appears to have a growth constant $s_+>1$.
Thus, for the convergence analysis presented in the next section, this norm sequence is too conservative to be useful.
On the other hand, the growth bounds above for the semi-norm sequence will provide the desired convergence results.

The next lemma is needed for the growth bounds for the sequence $T^{-1}_k$, $k\geq 0$.
In this instance, since all entries in each $T_k^{-1}$ are non-negative, $\snorm{T_k^{-1}}_{\infty}=\norm{T_k^{-1}}_{\infty}$ for all $k\geq 0$.

\begle \label{le:Tk-inv-growth}
For $\nu\in X^\ast$ and integers $i_1\geq 0$, \ldots, $i_n\geq 0$
such that $i_1+\cdots +i_n=\abs{\nu}$,
\begdi
\sum_{\xi_1\in X^{i_1},\ldots,\xi_n\in X^{i_n}} (\xi_1\shuffle\cdots\shuffle\xi_n, \nu) = \binom{\abs{\nu}}{i_1, \ldots, i_n},
\enddi
where the symbol on the right-hand side is the multinomial coefficient.
\endle

\begth \label{th:L-inv-growth}
The matrix representation $T_k^{-1}$ of the $\re$-linear map ${\mathscr L}^{-1}_k$ has
growth bounds
\begdi
k!\leq\norm{T_k^{-1}}_{\infty}\leq \left(\frac{1}{\ln 2}\right)^{k} k!,\;\;\forall k \geq 0.
\enddi
\endth

\begpr
The $\eta_i$-th row sum of $T^{-1}_k$ is by definition
\begin{align*}
\sum_{\Lbasis_j\in L^k} [T^{-1}]_{\eta_i,\Lbasis_j}&=\sum_{\Lbasis_j\in L^k} ({\mathscr L^{-1}}(\Lbasis_j),\eta_i),\;\;\eta_i\in X^k \\
&=\sum_{\Lbasis_j=l_{j_1}\cdots l_{j_n}\in L^k}( (l_{j_1})\shuffle\cdots\shuffle(l_{j_n}),\eta_i),\\
\end{align*}
where $\eta_i\in X^k$.
The lower bound on $\norm{T_k^{-1}}_{\infty}$ comes from the rows of $T_k^{-1}$ indexed by a word of the form $x_j^k$
with $(\Lbasis_j)=x_j$ so that
\begin{align*}
[T^{-1}]_{x_j^k,\ell_j^k} &= ({\mathscr L}^{-1}_k(\Lbasis_j^k),x_{j}^{k}) \\
&=((\Lbasis_j)^{\shuffle k}, x_j^{k}) = k!.
\end{align*}
The upper bound on $\norm{T_k^{-1}}_{\infty}$ comes from applying Lemma~\ref{le:Tk-inv-growth}, specifically,
\begin{align*}
\sum_{\Lbasis_j\in L^k} [T^{-1}]_{\eta_i,\Lbasis_j}
 &=\sum_{\Lbasis_j=l_{j_1}\cdots l_{j_n}\in L^k}( (l_{j_1})\shuffle\cdots\shuffle(l_{j_n}),\eta_i)\\
&\leq \sum_{n=1}^k
\sum_{i_1,\ldots,i_n> 0\atop k = i_1 + \cdots + i_n}\sum_{\xi_1\in X^{i_1},\ldots,\xi_n\in X^{i_n} \atop \abs{\eta_i} = i_1 + \cdots + i_n} (\xi_1\shuffle\cdots\shuffle\xi_n, \eta_i) \\
&=\sum_{n=1}^k \sum_{i_1,\ldots,i_n>0\atop k = i_1 + \cdots + i_n} \binom{k}{i_1, \ldots, i_n}.
\end{align*}
Note that the summation over the $i_v$'s above assumes $i_v>0$ since Lyndon words are never empty. In addition, $\xi_v\in X^{i_v}$ may not
always be a Lyndon word. This explains the inequality in the second line. The more delicate task is to evaluate the final summation above,
which cannot be done via the usual multinomial theorem.
Instead, observe that
\begin{align*}
\sum_{i_1,\ldots,i_n>0 \atop k = i_1 + \cdots + i_n} \binom{k}{i_1, \ldots, i_n}&=((\expup^z-1)^n,z^k) \\
&=n!\frac{((\expup^z-1)^n,z^k)}{n!} \\
&=n!\left\{k \atop n\right\},
\end{align*}
where $\left\{k \atop n\right\}$ denotes a Stirling number of the second kind. Therefore,
\begeq \label{eq:OBN-bound-for-invT}
\sum_{\Lbasis_j\in L^k} [T^{-1}]_{\eta_i,\Lbasis_j} \leq
\sum_{n=1}^k n!\left\{k \atop n\right\}.
\endeq
The integer sequence
\begdi
a(k)=\sum_{n=0}^k n!\left\{k \atop n\right\},\;\;k\geq 1
\enddi
is called the {\em ordered Bell numbers}. Since $\left\{k \atop n\right\}=0$ for $n=0$ and $k\geq 1$, the integer sequence
in \rref{eq:OBN-bound-for-invT} bounding the matrix norms is simply this sequence without the first term $a(0)=1$.
The ordered Bell numbers is sequence A000670 in the OEIS database. The first few terms are
\begin{align*}
\{a(k):k\geq 0\}&=\{1, 1, 3, 13, 75, 541, 4683, 47293,
545835, 7087261, 102247563,\ldots\}.
\end{align*}
The asymptotic behavior of this sequence is known to be
\begdi
a(k)\sim \frac{1}{2}\left(\frac{1}{\ln(2)}\right)^{k+1}k!,\;\;k\gg 1.
\enddi
From this result, it can be verified that
$\norm{T_k^{-1}}_{\infty}\leq \left(\frac{1}{\ln 2}\right)^{k} k!$, $k \geq 0$.
\endpr

The data suggests that $\norm{T_k^{-1}}_{\infty}=k!$, $k\geq 0$. But the limit given above has
the geometric growth constant that $1/\ln 2=1.4426950>1$. So it remains an open problem how to verify this
tighter upper bound.

\subsection{Local convergence under Lyndon word transduction}

The main result of this section is that
the Lyndon word transduction presented in Theorem~\ref{th:Lyndon-word-transduction}
preserves local convergence. But first it is necessary to establish sufficient conditions under which
a Chen--Fliess series indexed by Lyndon monomials is locally convergent. Consider the following
theorem.

\begth \label{th:LC-on-R[[L]]}
Suppose $c_L\in\allseriesL\cup\re\emptyset$ is a series with coefficients that satisfy
\begdi
|(c_L,l)|\leq K_L M_L^{\abs{l}},\;\;\forall l\in L^\ast
\enddi
for some real numbers $K_L,M_L>0$.
Then there exists real numbers $R_L,\Delta_L>0$ such that
for each $u\in B_1^m(R_L)[t_0, t_0+\Delta_L]$, the series
\begdi 
y(t) = F_{c_L}[u](t) = \sum_{l\in L^{\ast}}(c_L,l)\,E_{{\mathscr L}^{-1}(l)}[u](t,t_0)
\enddi
converges absolutely and uniformly on $[t_0, t_0+\Delta_L]$.
\endth

\begpr
Without loss of generality, assume $t_0=0$ and abbreviate $E_{\eta}(t,0)$ by $E_\eta(t)$.
For any $u\in L_1^m[0,\Delta_L]$ define $\bar{u}_j(t)=\abs{u_j(t)}$, $j=1,2,\ldots,m$ so that for all
$\eta\in X^\ast$
\begdi
\abs{E_\eta[u](t)}\leq E_\eta[\bar{u}](t),\;\;0\leq t\leq \Delta_L.
\enddi
It is known that for any $k\geq 0$ and $u\in L_1^m[0,\Delta_L]$ satisfying $\max\{\norm{u}_1,\Delta_L\}\leq R_L$ that
\begeq \label{eq:kfac-sumXk-E-bound}
k!\sum_{\eta\in X^k} \abs{E_\eta[u](t)}\leq (R_L(m+1))^k,\;\;t\in[0,\Delta_L]
\endeq
(see, for example, \cite{Duffaut_Espinosa-etal_09}).
Observe
\begin{align*}
\abs{F_{c_L}[u](t)}&\leq\sum_{l\in L^\ast} |(c_L,l)\,E_{{\mathscr L}^{-1}(l)}[u](t)| \\
&\leq \sum_{k=0}^\infty \sum_{l\in L^k} |(c_L,l)|\abs{E_{{\mathscr L}^{-1}(l)}[u](t)} \\
&\leq \sum_{k=0}^\infty \sum_{l\in L^k} K_LM_L^{|l|}\abs{E_{{\mathscr L}^{-1}(l)}[u](t)}.
\end{align*}
Applying \rref{eq:kfac-sumXk-E-bound} and Theorem~\ref{th:L-inv-growth}, it follows that
\begin{align*}
\abs{F_{c_L}[u](t)}&\leq K_L\sum_{k=0}^\infty M_L^k\sum_{l\in L^k}
\, \sum_{\eta\in X^k} ({\mathscr L}^{-1}(l),\eta) E_{\eta}[\bar{u}](t) \\
&\leq K_L\sum_{k=0}^\infty M_L^k
\sum_{\eta\in X^k}\left(\sum_{l\in L^k}
[T^{-1}_k]_{\eta,l}\right) E_{\eta}[\bar{u}](t) \\
&\leq K_L\sum_{k=0}^\infty M_L^k \norm{T^{-1}_k}_\infty
\sum_{\eta\in X^k} E_{\eta}[\bar{u}](t) \\
&\leq \frac{K_L}{\sqrt{2\pi}}\sum_{k=0}^\infty (\expup M_LR_L(m+1))^k.
\end{align*}
Thus, the series converges provided $\expup M_LR_L(m+1)<1$, and the claims in the theorem then follow directly.
\endpr

\begex \label{ex:CFS-in-Lyndon-coordinates}
Suppose $X=\{x_0,x_1\}$ and $M_L>0$ is given.
If $c_L=\sum_{k\geq 0} M_L^k l_1^k$, then
\begin{align*}
F_{c_L}[u]&=\sum_{k=0}^\infty M_L^k E_{{\mathscr L}^{-1}(l_1^k)}[u] \\
&=\sum_{k=0}^\infty M_L^k E_{(l_1)}^k[u] \\
&=\sum_{k=0}^\infty M_L^k E_{x_1^{\shuffle\,k}}[u] \\
&=\sum_{k=0}^\infty M_L^k k!\,E_{x_1^k}[u].
\end{align*}
Therefore, on any interval $[0,\Delta_L]$
\begin{align*}
\abs{F_{c_L}[u](t)}&\leq \sum_{k=0}^\infty M_L^k k!\frac{\norm{u}_1^k}{k!} \\
&=\frac{1}{1-M_L\norm{u}}_1
\end{align*}
when $\norm{u}_1<1/M_L$.
\endex

Now the main convergence preservation theorem is presented.

\begth
For every $c\in\allseriesLC$, there exists a $\Delta>0$ and a ball $B_1^m(R)[0,\Delta]$ of positive radius $R$ in $L_1^m[0,\Delta]$ such that
$
F_{{\mathscr L}(c)}[u]
$
converges absolutely and uniformly on $[0,\Delta]$
for all $u\in B_1^m(R)[0,\Delta]$.
\endth

\begpr
Given that $c$ satisfies the growth rate \rref{eq:local-convergence-growth-bound},
it follows from Theorem~\ref{th:L-growth} that for any $l\in L^\ast$
\begin{align*}
\abs{({\mathscr{L}(c)},l)}&=\abs{\sum_{\eta\in X^\ast} (c,\eta)({\mathscr L}(\eta),l)}\\
&\leq K M^{\abs{l}} \abs{l}! \abs{\sum_{\eta\in X^{\abs{l}}} ({\mathscr L}(\eta),l)} \\
&=K M^{\abs{l}} \abs{l}! \abs{\sum_{\eta\in X^{\abs{l}}} [T_{\abs{l}}]_{l,\eta}} \\
&\leq K M^{\abs{l}} \abs{l}!\,\snorm{T_{\abs{l}}}_{\infty} \\
&\leq \frac{K}{\sqrt{\pi/2}} (2M)^{\abs{l}}.
\end{align*}
The claim is proved in light of Theorem~\ref{th:LC-on-R[[L]]}.
\endpr

\begex
Consider the locally convergent generating series $c=\sum_{k\geq 0} M^kk!\,x_1^k$.
Then $F_c[u]=F_{{\mathscr L}(c)}[u]$,
where
\begin{align*}
{\mathscr L}(c) &= \sum_{k=0}^\infty M^k k!\, {\mathscr L}(x_1^k) \\
&=\sum_{k=0}^\infty M^k k! \frac{l_1^k}{k!} \\
&=\sum_{k=0}^\infty M^k l_1^k,
\end{align*}
which is the generating series in Example~\ref{ex:CFS-in-Lyndon-coordinates} with $M_L=M$. So
the convergence analysis presented there applies here.
\endex

\begex
Consider the locally convergent generating series $c=\sum_{\eta\in X} KM^{\abs{\eta}}\abs{\eta}!\,\eta$.
From Lemma~\ref{le:Tk-growth}, it follows that
\begin{align*}
F_{{\mathscr L}(c)}[u]&=\sum_{k=0}^\infty KM^k F_{{\mathscr L}(\charseries(X^k))}[u] \\
&=K\sum_{k=0}^\infty (M F_{\charseries(L^1)}[u])^k \\
&=\frac{K}{1-M F_{\charseries(L^1)}[u]},
\end{align*}
which converges
for sufficiently small $u$.
\endex

\section{Improving the Efficiency of Computing Chen--Fliess Series}
\label{sec:computational-efficiency}

In this section, the concept of computational efficiency is introduced in order
to determine a priori how beneficial a change of basis for a Fliess operator evaluation could be from
the point of view of reducing the number of integrals that need to be computed.
It is not meant to apply to a particular Fliess operator, nor is it
a precise indicator of the computational complexity of a specific integration algorithm.
Given this generality, it is nevertheless not trivial to compute
exactly. But it will be shown that reasonable estimates of this performance
parameter are possible, and they provide some guidance
in this regard. The concept will be demonstrated in the
context of a cyber-physical security application
involving a continuous stirred-tank reactor system.

\subsection{Computational efficiency}

For $p\in\allpoly$, let $J_X(p)$ be the minimum number of integrals required to compute $F_p$ in terms of iterated integrals indexed by words over $X$.
No existing algebraic relationships between iterated integrals are used except for the fact that integrals for longer words rely on integrals for shorter words.
For example, if $p_1=x_i^2$ and $p_2=x_i+x_i^2$, then $J_X(p_1)=J_X(p_2)=2$.
The worst case scenario for polynomials
of degree $n\geq 1$ in $\allpoly$ is denoted by $I_X(n)=\max_{\deg(p)=n} J_X(p)$.

\begth
Let $X$ be a finite alphabet. If $p\in\allpoly$ with $\deg(p)\geq 1$,
then
\begdi
\deg(p)\leq J_X(p)\leq I_X(\deg(p)),
\enddi
where
\begdi
I_X(n)=\sum_{k=1}^n (\card(X))^k.
\enddi
In particular, if $\card(X)=1$, then $J_X(p)=\deg(p)$ so that $I_X(n)=n$, $n\geq 1$. Otherwise,
\begeq \label{eq:IXplus}
I_X(n)=\frac{\card(X)((\card(X))^n-1)}{\card(X)-1},\;\;n\geq 1.
\endeq
\endth

(Note this upper bound appears in \cite[p.7]{Reizenstein_17}.)

\begpr
If $\card(X)=1$ and $\deg(p)=n\geq 1$, then necessarily $p=\sum_{k=1}^n (p,x_i^k)x_i^k$, where $X=\{x_i\}$ and $(p,x_i^n)\neq 0$.
Thus, $J_X(p)=\deg(p)$.
If $\card(X)>1$, then the lower bound on $J_X(p)$ is established in an analogous fashion.
If $p=\charseries(X^n)$, then  $J_X(p) = \sum_{k=1}^n(\card(X))^{k}$. All other cases
fall between these two extremes.
The expression for $I_X(n)$ in this case follows from a standard formula for geometric series.
\endpr

Let $q\in\re[L]\cup\re\emptyset$ and $J_L(q)$ be the minimum number of integrals required to
compute $F_{q}$ in terms of iterated integrals indexed by Lyndon words. Again, the only assumption made in the counting process is that
integrals for longer words will rely on integrals for shorter words.
The worst case scenario for polynomials
of degree $n\geq 1$ is denoted by $I_L(n)=\max_{\deg(q)=n} J_L(q)$.
For example, if $\card(X)=1$, then trivially, $J_L(q)=1$ so that $I_L(n)=1$, $n\geq 1$.
For larger alphabets, computing $I_L(n)$ is more difficult, as described in the next example.

\begin{table}[tb]
\renewcommand{\arraystretch}{1.5}
\begin{center}
\caption{Computing $I_L(n)$ (red) in  Example~\ref{ex:computing-ILplus-n}}
\label{tbl:nest-L-integral-degree-3}
{\scriptsize
\begin{tabular}{|c|c|c|c|} \hline
$n$ & $l_i$ (previous integration) & number of new integrations  & running sum \\ \hline\hline
\multirow{2}[0]*{1} & $x_0$ & 1 & 1 \\ \cline{2-4}
& $x_1$ & 1 & \textcolor[rgb]{1.00,0.00,0.00}{2} \\ \hline\hline
2& $x_0(x_1)$ & 1 & \textcolor[rgb]{1.00,0.00,0.00}{3} \\ \hline\hline
\multirow{2}[0]*{3}& $x_0(x_0x_1)$ & 1 & 4 \\ \cline{2-4}
& $x_0x_1(x_1)$ & 2 & \textcolor[rgb]{1.00,0.00,0.00}{6} \\ \hline\hline
\multirow{3}[0]*{4}& $x_0 (x_0 x_0 x_1)$      &    1 & 7 \\ \cline{2-4}
& $x_0 (x_0 x_1 x_1)$      &    1 & 8 \\ \cline{2-4}
& $x_0 x_1 (x_1 x_1)$      &    2 & \textcolor[rgb]{1.00,0.00,0.00}{10} \\  \hline\hline
\multirow{6}[0]*{5}& $x_0 (x_0 x_0 x_0 x_1)$  & 	1 &  11 \\  \cline{2-4}
& $x_0 (x_0 x_0 x_1 x_1)$  & 	1 & 12 \\ \cline{2-4}
& $x_0 x_0 x_1 (x_0 x_1)$  & 	3 & 15 \\ \cline{2-4}
& $x_0 (x_0 x_1 x_1 x_1)$  & 	1 & 16 \\ \cline{2-4}
& $x_0 x_1 (x_0 x_1 x_1)$  & 	2 & 18 \\ \cline{2-4}
& $x_0 x_1 (x_1 x_1 x_1)$  & 	2 & \textcolor[rgb]{1.00,0.00,0.00}{20} \\ \hline
\end{tabular}
}
\end{center}
\end{table}

\begex \label{ex:computing-ILplus-n}
Suppose $X=\{x_0,x_1\}$ and $q_n=\sum_{k=1}^n\charseries(L^k)$.
Then the iterated integral for every Lyndon word up to length $k$ will need to be computed.
Table~\ref{tbl:nest-L-integral-degree-3} shows
how to compute
$I_L(n)=J_L(q_n)$ up to length $n=5$ (shown in red). Each time the length of a word increases by one letter, at least one additional integral needs to be
computed. In five of the 14 cases, more than one integral is required. For example, the word $x_0x_1x_1x_1$ requires two additional integrals because $E_{x_1x_1}$
was computed earlier for the word $x_0x_1x_1$. Note, in particular, that $x_1x_1$ is {\em not} a Lyndon word,
so $x_0x_1(x_1x_1)$ is not a factorization in terms of Lyndon words.
In fact, the Lyndon word factorization of $x_0x_1x_1x_1$ with the longest right factor, i.e., the {\em standard factorization}, is $x_0x_1x_1(x_1)$.
Assuming that previous integrals are available only for shorter Lyndon words, this word would require $3>2$ additional integrals.
So it is not obvious whether any algorithm based solely on Lyndon word factorizations
can render the integer sequence $I_L(n)$, $n\geq 1$.
A brute force counting approach, as illustrated above, yields
\begin{align*}
\{I_L(n): n\geq 1\}&=\{
2,3,6,10,20,33,66,116,222,406,778,1442,2762,5203,9952,\cdots \\
&\hspace{0.3in} 18941,36348,69575,133930,257646
\}.
\end{align*}
This sequence is not in the OEIS database at present, and to the authors' best knowledge,
an analytical formula for computing this sequence is not known.
\endex

\begin{figure}[tb]
\begin{center}
\includegraphics[scale=0.55]{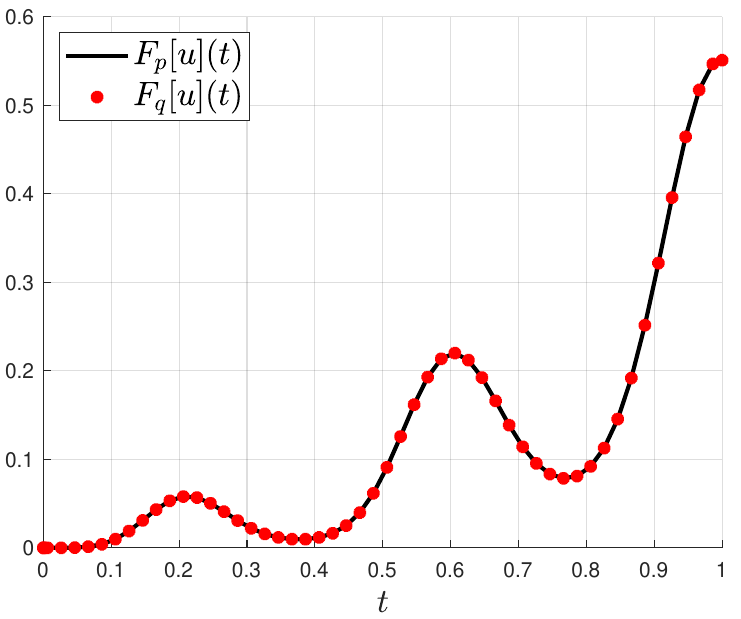}
\caption{$F_p[u]$ and $F_{q}[u]$ in Example~\ref{ex:Lyndon-poly-simulation}, where $p=\charseries(X^3)$, $q={\mathscr L}(p)$, and $u(t)=4\sin(16t)$}
\label{fig:Lambert-series-simulation}
\end{center}
\end{figure}

\begex \label{ex:Lyndon-poly-simulation}
Let $X=\{x_0,x_1\}$, $p=\charseries(X^3)$, and $q={\mathscr L}(p)=(l_0+l_1)^3/3!$ using Lemma~\ref{le:Tk-growth}.
Then $\deg(p)=3$, $J_X(p)=I_X(3)=14$, $J_L(q)=2$, and $I_L(3)=6$. Therefore,
$\deg(p)<J_X(p)=I_X(3)$, $J_L(q)<I_L(3)$, and
$J_L(q)<J_X(p)$.
A MATLAB simulation of both $F_p[u]$ and $F_q[u]$ when $u(t)=4\sin(16t)$
is shown in Figure~\ref{fig:Lambert-series-simulation}. As expected, the outputs match exactly.
\endex

\begex
In applications, where usually $\deg(p)\gg 1$, it is common to have
$J_L({\mathscr L}(p))\ll J_X(p)$.
There are, however, exceptional cases where $J_X(p)< J_L({\mathscr L}(p))$.
For example, if $p=x_1x_0$, then $q={\mathscr L}(p)=-l_2+l_1 l_0$ so that
$J_X(p)=2<J_L(q)=I_L(2)=3$.
\endex

For any $n\geq 1$, the quantity
\begdi
\compeff(n):=\frac{I_X(n)-I_L(n)}{I_X(n)}
\enddi
will be called the
{\em computational efficiency} for degree $n$.
It is clearly the relative difference between the upper bounds on $J_X(p)$ and $J_L({\mathscr L}(p))$ when $\deg(p)=n$. It will be viewed as the
maximum possible relative reduction in the number of integrals needed to compute $F_p[u]$
as a result of applying the Chen--Fliess transduction ${\mathscr L}^\dag$.
For example, when $\card(X)=1$, it follows that $\compeff(n)=1-(1/n)$ for all $n\geq 1$.
For the general case, the following theorem provides bounds for $\compeff(n)$.

\begth \label{th:bounds-compeff}
If $\card(X)\geq 2$ and finite, then
\begdi
\compeff_-(n)\leq \compeff(n)\leq \compeff_+(n),\;\;n\geq 1,
\enddi
where
\begin{align*}
\compeff_-(n)&:=\frac{I_X(n)-(nL(n)+1)}{I_X(n)} \\
\compeff_+(n)&:=\frac{I_X(n)-L_+(n)}{I_X(n)}.
\end{align*}
\endth

\begpr
The claim follows from the assertion that
\begdi
L_+(n)\leq I_L(n)\leq nL(n)+1,\;\; n\geq 1.
\enddi
The lower bound follows from the conservative assumption that as the word length increases from $n$ to $n+1$, only a single new integral needs to be computed for
each additional Lyndon word.
As shown in Example~\ref{ex:computing-ILplus-n}, this happens often but not always. So in general, $L_+(n)$ will under count
the number of required integrals as a function of $n$.
To derive the upper bound, a suboptimal scheme based on Lyndon word factorizations is used.
A well known alternative characterization of Lyndon words is that
a word $\eta\in X^+$ of length greater than one is a Lyndon word if and only if $\eta=l_il_j$ with $l_i,l_j\in L$ and $l_i<l_j$
\cite[Proposition 5.1.3]{Lothaire_83}. Thus, given a Lyndon word of length $n>1$, every Lyndon word $l_j$ will appear as a right factor of
a Lyndon word in $L(n)$ with one exception, the smallest Lyndon word $l_0$ will never appear as a right factor.
Thus, it is necessary to only consider the integrals for words in $L(n)$, each of which requires $n$ integrals, and
then to add one to compensate for the missing integral $E_{x_0}$ when $n>1$. This results in the given upper bound.
\endpr

The following is an asymptotic version of Theorem~\ref{th:bounds-compeff}.

\begth \label{th:asymptotic-bounds-compeff}
If $\card(X)\geq 2$, then for all $n\gg 1$,
\begdi
0< \hatcompeff_-(n)\leq \compeff(n) \leq  \hatcompeff_+(n)<1,
\enddi
where
\begin{align*}
\hatcompeff_-(n)&:=\frac{1}{\card(X)} \\
\hatcompeff_+(n)&:=1-\frac{1}{n}(\card(X)-1).
\end{align*}
\endth

\begpr
The upper bound is proved first. From \rref{eq:IXplus} and the asymptotic approximation in
Lemma~\ref{le:Lyndon-word-asymptotics}, it follows directly that if $n\gg 1$, then
\begin{align*}
\compeff_+(n)&\sim 1-\frac{1}{n} \frac{(\card(X))^n}{(\card(X))^n-1}(\card(X)-1)\\
&\sim 1-\frac{1}{n}(\card(X)-1) \\
&=:\hatcompeff_+(n).
\end{align*}
Clearly, if $n\gg 1$, then $\hatcompeff_+(n)<1$. The lower bounds follow in a similar fashion.
\endpr

\begin{figure}[tb]
\begin{center}
\includegraphics[scale=0.55]{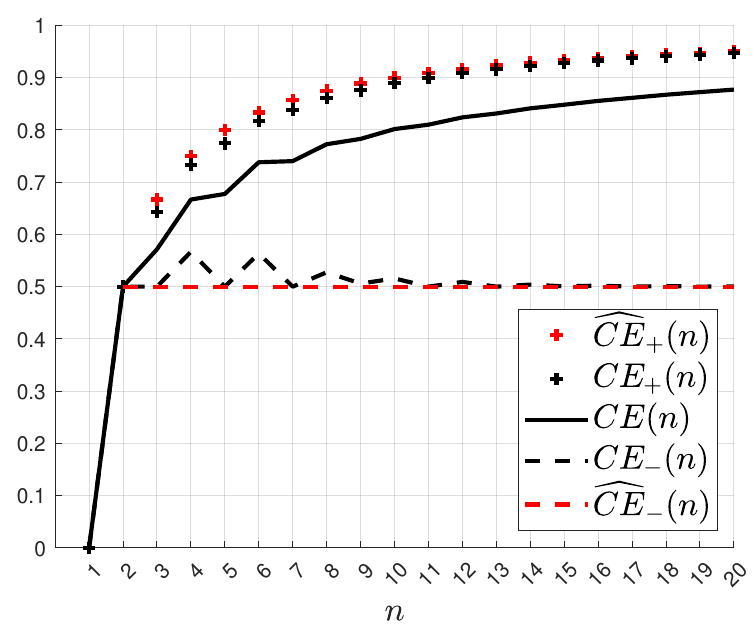}
\caption{$\compeff(n)$ and its various bounds when $\card(X)=2$}
\label{fig:compeff-bounds}
\end{center}
\end{figure}

\begin{figure}[tb]
\begin{center}
\includegraphics[scale=0.55]{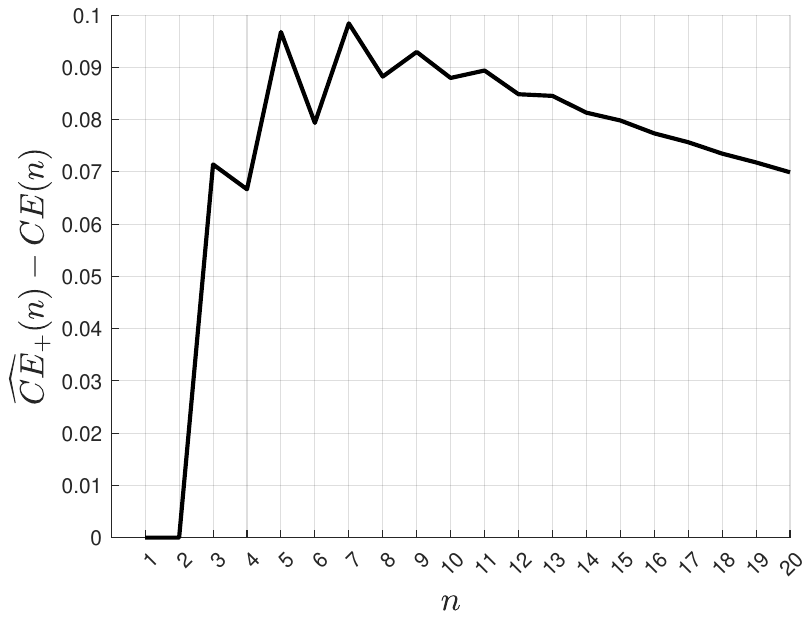}
\caption{$\hatcompeff_+(n)-\compeff(n)$ when $\card(X)=2$}
\label{fig:Ceff-error}
\end{center}
\end{figure}

\begin{figure}[tb]
\begin{center}
\includegraphics[scale=0.55]{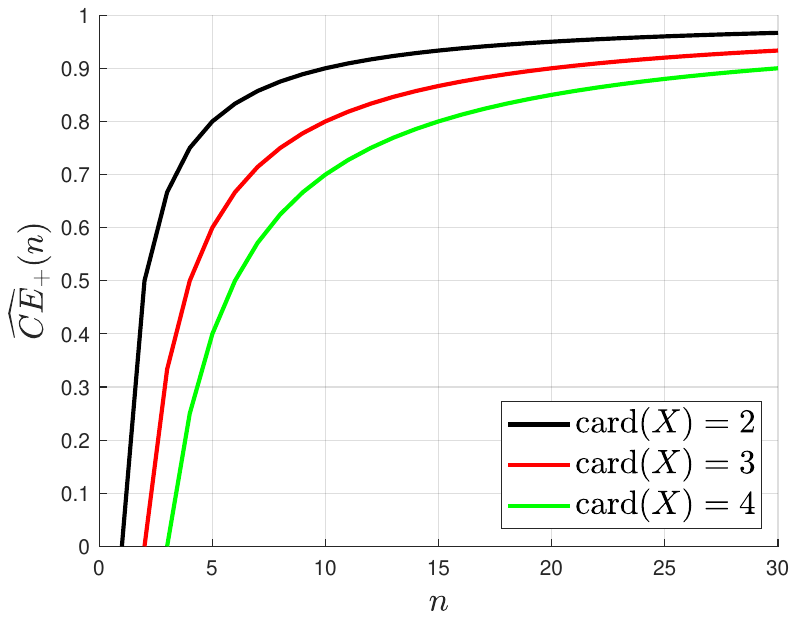}
\caption{$\hatcompeff_+(n)$ when $\card(X)=2,3,4$}
\label{fig:Ec-varying-cardX}
\end{center}
\end{figure}

\begex
For the case where $\card(X)=2$, the computational efficiency and its various bounds as a function of $n$ are shown in
Figure~\ref{fig:compeff-bounds}. The lower bound is very conservative, but the upper bound has less than a 10\% error
as shown in Figure~\ref{fig:Ceff-error}.
So it is an easy to compute estimate of the
computational efficiency.
In this particular case, it indicates that for $n$ on the order of ten or higher, the number of integrals that need to be computed could
be reduced by as much as 90 percent.
Figure~\ref{fig:Ec-varying-cardX} shows this upper bound when $\card(X)=2,3,4$.
The true computational efficiency likely follows a similar trend.
As $\card(X)$ increases, there
are many more words in a polynomial of a given degree.
So initially $\compeff(n)$
is significantly lower as $\card(X)$ increases until $n$ becomes sufficiently large.
\endex

\subsection{Zero dynamics attack on a continuously stirred-tank reactor system}

Consider a first-order, exothermic, irreversible reaction of a reactant in a product substance
carried out in a well mixed continuously stirred-tank reactor (CSTR) system.
The mass and energy balances give
the dynamics (in dimensionless form):
\begin{subequations} \label{eq:CSRT_dynamics}
\begin{align}
\dot{z}&= \left[\begin{array}{c}
-z_1 + \alpha(1-z_1)\expup^{\frac{z_2}{1+z_2/\gamma}} \\
-(\beta+1)z_2 + \kappa \alpha(1-z_1)\expup^{\frac{z_2}{1+z_2/\gamma}}
\end{array}\right] +
\left[\begin{array}{c}
0 \\
\beta
\end{array}\right] u \\
y &= z_2.
\end{align}
\end{subequations}
Here $z_1$ is the reactant concentration,
$z_2$ is the reactor temperature, and
$u$ is the cooling reactor jacket temperature
\cite{Doyle-Henson_97,Uppal-etal_74}.
The physical constants $\alpha$, $\beta$, $\gamma$, and $\kappa$ are all
set to unity for convenience. For $z(0)=0$,
the corresponding generating series over $X=\{x_0,x_1\}$ for the input-output map $y=F_c[u]$ is computed directly from \rref{eq:CSRT_dynamics} and
found to be
\begin{align*}
c &=  x_0 + x_1 - 2 x_0^2 - x_0 x_1 - 2 x_0^2x_1 - 2 x_0x_1x_0
- x_0x_1^2 + 22 x_0^4 + 15 x_0^3x_1 \\
&\hspace*{0.12in}
+ 11 x_0^2x_1 x_0 + 4 x_0^2x_1^2
+ 6 x_0 x_1 x_0^2 + 2x_0x_1x_0x_1 + 2 x_0x_1^2x_0 + x_0x_1^3 + \cdots.
\end{align*}
In which case,
\begin{align*}
{\mathscr L}(c)&=l_0+l_1-l_0^2-l_2 - 2 l_0 l_2 + 2l_3 - l_4  +\frac{11}{12}l_0^4 + 3 l_0^2l_2
 - l_0l_3 + 2 l_0l_4  + l_7  - \frac{13}{15}l_0^5\\
&\hspace*{0.12in}
 - \frac{7}{3}l_0^3l_2 + \frac{3}{2}l_0l_2^2 +\frac{11}{2} l_0^2l_3
 - l_2l_3 - 2 l_0^2l_4  - 10 l_0l_5 +  3 l_0l_6 + 3l_8 - 2 l_9 + \cdots.
\end{align*}
In the context of cyber-physical security, a {\em zero-dynamics attack} is one where an attack input $u^\ast$ is injected into the system
to alter its internal (unmonitored) states while leaving the measured output $y$ unchanged \cite{Cao-etal_20,Jeon-Eun_19,Park-etal_18}.
If the target system is initially operating
in equilibrium, then without loss of generality this output can be assumed to be the zero function.
It was shown in \cite{Gray-etal_CISS21} that the zero dynamics attack input for this normalized CSTR system is
$u^\ast(t)=-(1+\expup^{-2t})/2$, $t\geq 0$. If the attacker has no knowledge of the system's dynamics \rref{eq:CSRT_dynamics}, its
generating series $c$ can be estimated using input-output data \cite{Gray-etal_ACC22}.
The objective here is to validate that $F_c[u^\ast]=0$ by
computing the Chen--Fliess series over some finite interval, where $c$ is suitably truncated.

To determine whether this CSTR system is a good candidate for applying
the Lyndon word transduction ${\mathscr L}^\dagger$, the computational efficiency bound for $n=9$
was computed for the two letter case, namely, $\hatcompeff_+(9)=8/9\approx 0.8889$, which represents a considerable potential savings.
To see what the exact computational consequences are of changing the Chen--Fliess series basis,
the system was simulated using Python 3.10 implementations of Algorithms~\ref{alg:algorithm1} and \ref{alg:algorithm2}.

\begin{algorithm}[h]
		\caption{Compute $y(t)=F_c[u](t)$ on $[t_0,t_1]$}
		\label{alg:algorithm1}
		\begin{algorithmic}[1]
			\renewcommand{\algorithmicrequire}{\textbf{Input:}}
			\renewcommand{\algorithmicensure}{\textbf{Output:}}
\REQUIRE $c$, $n$, $u$, $[t_0,t_1]$
\STATE Truncate $c$ to length $n$.
\STATE Compute  $E_{\eta}[u](t,t_0)$, $\forall t\in[t_0,t_1]$, $\eta\in X^{\le n}$.
\STATE Compute $y(t)=\sum_{\eta\in X^{\leq n}} (c,\eta) E_\eta[u](t,t_0)$, $\forall t\in[t_0,t_1]$.
\ENSURE $y(t)$, $t\in[t_0,t_1]$
		\end{algorithmic}
\end{algorithm}

\begin{algorithm}[h]
		\caption{Compute $y(t)= F_{{\mathscr L}(c)}[u](t)$ on $[t_0,t_1]$}
		\label{alg:algorithm2}
		\begin{algorithmic}[1]
			\renewcommand{\algorithmicrequire}{\textbf{Input:}}
			\renewcommand{\algorithmicensure}{\textbf{Output:}}
\REQUIRE $c$, $n$, $T$, $u$, $[t_0,t_1]$
\STATE Truncate $c$ to length $n$.
\STATE Compute coefficients $(\mathscr L(c),l)$ for $l=l_{i_1}\cdots l_{i_k}\in L^{\le n}$ using matrix representation $T$.
\STATE Compute  $E_{(l_i)}[u](t,t_0)$, $\forall t\in[t_0,t_1]$, $l_i\in L^{\le n}$.
\STATE Compute $E_{{\mathscr L}^{-1}(l)}[u](t,t_0)$ 
  in terms of $E_{(l_i)}[u](t,t_0)$.
\STATE Compute $y(t)=\sum_{l\in L^{\leq n}} (\mathscr L(c),l)E_{{\mathscr L}^{-1}(l)}[u] 
(t,t_0)$, $\forall t\in[t_0,t_1]$.
\ENSURE $y(t)$, $t\in[t_0,t_1]$
		\end{algorithmic}
\end{algorithm}

\begin{figure}[tb]
\hspace*{0.5cm}\includegraphics[scale=0.5]{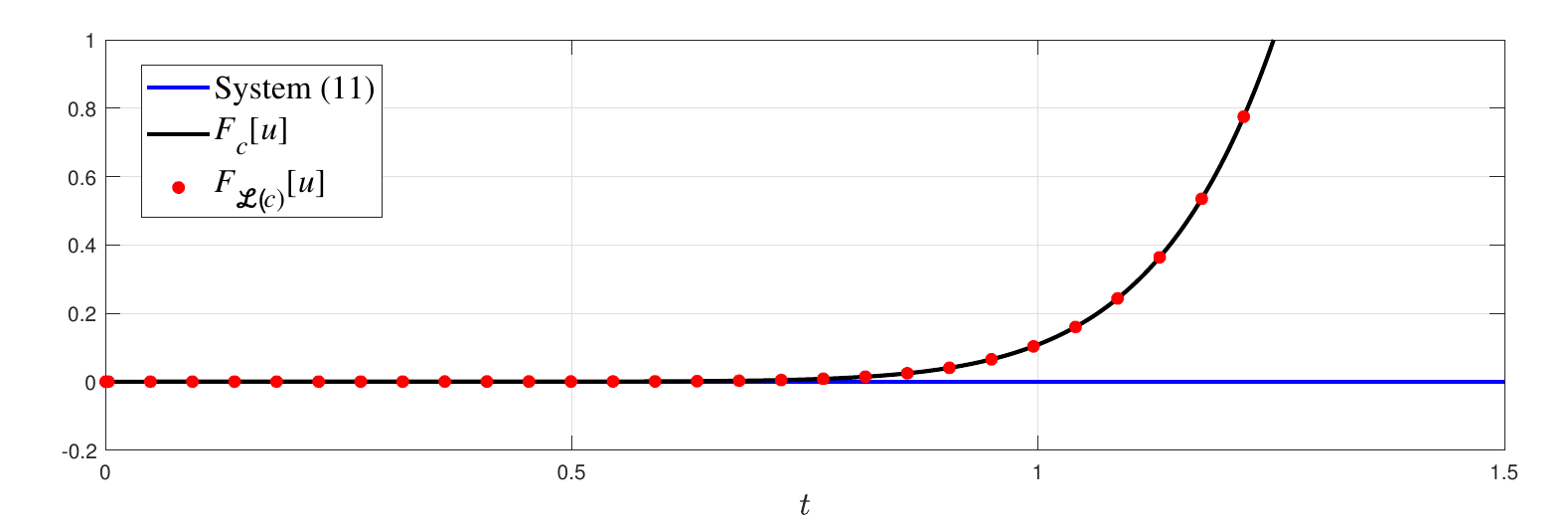}
\caption{Simulated responses of the CSTR system with $n=9$ for the
Chen--Fliess series}
\label{fig:CSTR-simulations}
\end{figure}

For each algorithm, the method for computing an iterated integral $E_{\eta}[u](t,t_0)$, $\eta\in X^\ast$ employed \emph{Chen's identity} \cite{Chen_57,Chen_61} and only first-order path increments. Specifically, on each subinterval $[t_j,t_{j+1}]$, the integral indexed by a word $\eta = x_{i_1}\cdots x_{i_k}$ reduces to a scaled product of increments,
\begdi  
E_{\eta}[u](t_{j+1},t_j) \;=\; \frac{1}{k!}\prod_{r=1}^k \Delta U^{(i_r)}_{t_j,t_{j+1}},
\enddi
where
\begin{equation*}
\Delta U^{(i)}_{t_j,t_{j+1}} \;=\; \int_{t_j}^{t_{j+1}} u_i(\tau)\,d\tau.
\end{equation*}
Iterating Chen’s identity across $M$ increments yields a structured polynomial expansion that replaces the
need for explicitly computing higher-order integrals of $u$ with arithmetic operations. For example, computing the iterated integral of a Lyndon word requires $\binom{k+M-1}{M-1}$ additions for words of length $k$ over a time interval with $M$ subintervals \cite{Reizenstein-Graham_20}.
Nevertheless, the more iterated integrals that need to be computed,
the heavier the computational load in terms of arithmetic operations.

The code was executed on an
Intel Core i7-7820HQ processor with base clock speed 2.9 GHz/turbo boost speed 3.9 GHz, quad-core, and 8 MB of cache memory.
The output was computed when $u=u^\ast$ using both algorithms for various maximum word lengths $n$.
The results are shown in Figure~\ref{fig:CSTR-simulations} when $n=9$. As expected, $F_c[u]=F_{{\mathscr L}(c)}[u]\approx 0$ over a finite interval, but the consequences of truncating the
generating series become evident
starting at $t=0.75$
when compared against a simulation of the state space realization \rref{eq:CSRT_dynamics}.
The simulation times for both algorithms were
measured using the {\tt time.time} command in the Python library {\tt time}. For Algorithm~\ref{alg:algorithm1}, steps $2$ to $3$ were timed, while for
Algorithm~\ref{alg:algorithm2}, steps $2$ to $5$ were timed. Memory usage  was computed using the {\tt asizeof} command in the Python library {\tt pympler}.
It returns the total memory size (in bytes) of a Python object, including the size of all objects it refers to.
The mean values of these performance parameters for 10 runs as a function of $n$ are shown in Table~\ref{tbl:CSRT-sim-statistics}.
Observe that for the $n=9$ case, the simulation time was reduced on average using Algorithm~\ref{alg:algorithm2} by a factor of $3.01$ as compared to
Algorithm~\ref{alg:algorithm1}. Similarly, the
memory usage was reduced by a factor of 4.12.

\begin{table}[tb]
\renewcommand{\arraystretch}{1.5}
\begin{center}
\caption{Mean simulation times and memory usages for the CSTR system models $F_c$ (Algorithm~\ref{alg:algorithm1})
and $F_{{\mathscr L}(c)}[u]$ (Algorithm~\ref{alg:algorithm2})  as a function of $n$ for 10 runs}
\label{tbl:CSRT-sim-statistics}
{\scriptsize
\begin{tabular}{|c||c|c||c|c|} \hline
max word &
$F_c$ mean &
$F_{{\mathscr L}(c)}[u]$ mean &
$F_c$ mean &
$F_{{\mathscr L}(c)}[u]$ mean \\[-0.03in]
length $n$ &
simulation time (sec) &
simulation time (sec) &
memory usage (MB) &
memory usage (MB) \\ \hline\hline
2 & 0.013 & 0.014 & 0.292 & 0.196 \\
\hline
3 & 0.032 & 0.027 & 0.560 & 0.260 \\
\hline
4 & 0.105 & 0.049 & 1.084 & 0.436 \\
\hline
5 & 0.208 & 0.088 & 2.145 & 0.800 \\
\hline
6 & 0.433 & 0.194 & 4.382 & 1.467 \\
\hline
7 & 0.940 & 0.316 & 8.741 & 2.267 \\
\hline
8 & 2.454 & 0.875 & 17.451 & 4.318 \\
\hline
9 & 4.868 & 1.616 & 34.885 & 8.472 \\
\hline
\end{tabular}
}
\end{center}
\end{table}

\section{Conclusions}

A computationally efficient way to evaluate Chen--Fliess series was presented by introducing a change of
basis for the computation using Lyndon words.
Specifically, it was shown that one can improve the efficiency of computing the output function
using the fact that
the shuffle algebra on (proper) polynomials over a finite alphabet is isomorphic to the
polynomial algebra generated by the Lyndon words over this alphabet.
Thus, the iterated integrals indexed by Lyndon monomials contain all the input information needed to compute the output.
The change of basis was accomplished by applying a transduction
to re-index the Chen--Fliess series in terms of Lyndon monomials.
It preserves the local
convergence property of the original series.
The method was illustrated by simulating a zero dynamics attack on a continuously stirred reactor tank system.

\section*{Acknowledgements}

L. A. Duffaut Espinosa was supported by
NSF CAREER Award 2340089 and
NSF EPSCoR Research Fellows Award 2429505.

\appendix

\section{Proof of Lemma~\ref{le:Tk-growth}}
\label{proof-of-lemma-le:Tk-growth}

From the identity
\begdi
\charseries(X^k)
= \sum_{r_0,r_1,\ldots,r_m\geq 0 \atop r_0+r_1+\cdots+r_m=k}
x_0^{r_0}\shuffle x_1^{r_1}\shuffle \cdots \shuffle x_m^{r_m}
\enddi
(see \cite{Duffaut_Espinosa-etal_09}), the multinomial theorem, and the identity $x_i^r=x_i^{\shuffle r}/r!$,
it follows that
\begin{align*}
{\mathscr L}(\charseries(X^k))&= \sum_{r_0,r_1,\ldots,r_m\geq 0 \atop r_0+r_1+\cdots+r_m=k}
{\mathscr L}(x_0^{r_0}) {\mathscr L}(x_1^{r_1})\cdots {\mathscr L}(x_m^{r_m}) \\
&= \sum_{r_0,r_1,\ldots,r_m\geq 0 \atop r_0+r_1+\cdots+r_m=k}
 \frac{l_0^{r_0}}{r_0!} \frac{l_1^{r_1}}{r_1!}\cdots  \frac{l_m^{r_m}}{r_m!} \\
&=(l_0+l_1+\cdots+l_m)^k\frac{1}{k!} \\
&=(\charseries({\mathscr L}(X)))^k\frac{1}{k!}.
\end{align*}

\section{Proof of Lemma~\ref{le:Tk-inv-growth}}
\label{proof-of-lemma-Tk-inv-growth}

The proof is by induction on the length of $\nu$.
The case $\abs\nu=0$, $i_1=0=\cdots=i_n=0$ is trivial.  For the case $\abs\nu=1$, let $i_1=1$ and $i_2=0=\cdots=i_n=0$ to get
\begdi
\sum_{\substack{{\xi_1\in X}\\{\xi_2,\ldots,\xi_n\in X^{0}}}} (\xi_1\shuffle\emptyset\shuffle\cdots\shuffle\emptyset, \nu) = 1 = \binom{1}{1,0,\ldots,0}.
\enddi
All other such cases are similar.
Thus, the identity holds for $\abs\nu=1$.
Now suppose the claim is true for $\nu^\prime\in X^{j-1}$ with $j\geq 2$, and let
$\nu=x\nu^{\prime}$ such that $\nu\in X^{j}$ with $x\in X$. Assume $i_1+\cdots + i_n=j$ for fixed integers $i_1\geq0,\ldots,i_n\geq0$.
The identity in question can always be rewritten without loss of generality in the form
\begdi
\sum_{\substack{{\xi_1\in X^{i_1}},{\ldots},{\xi_{n^\prime}\in X^{i_{n^\prime}}}}} (\xi_1\shuffle\cdots\shuffle\xi_{n^\prime}, \nu) = \binom{\abs\nu}{i_1,\ldots,i_{n^\prime}}
\enddi
where $i_1\geq 1,\ldots,i_{n^\prime}\geq 1$ and $i_{{n^\prime}+1}=\cdots=i_n=0$.
Observe
\begin{align*}
\sum_{\substack{{\xi_1\in X^{i_1}},{\ldots},{\xi_{n^\prime}\in X^{i_{n^\prime}}}}} (\xi_1\shuffle\cdots\shuffle\xi_{n^\prime}, \nu)
&= \sum_{\substack{{\xi_1\in X^{i_1}},{\ldots},{\xi_{n^\prime}\in X^{i_{n^\prime}}}}} (\xi_1\shuffle\cdots\shuffle\xi_{n^\prime}, x\nu^{\prime}) \\
&= \sum_{\substack{{\xi_1\in X^{i_1}},{\ldots},{\xi_{n^\prime}\in X^{i_{n^\prime}}}}} (x^{-1}(\xi_1\shuffle\cdots\shuffle\xi_{n^\prime}),\nu^{\prime}) \\
&=\sum_{l=1}^{n^{\prime}} \sum_{\substack{{x^{\prime}\in X},\\{\xi_1\in X^{i_1}},{\ldots},\\{\xi_l^{\prime}\in X^{i_l-1}},{\ldots},\\{\xi_{n^\prime}\in X^{i_{n^\prime}}}}}
(\xi_1\shuffle\cdots\shuffle x^{-1}(x^{\prime}\xi_l^{\prime})\shuffle\cdots\shuffle\xi_{n^\prime},\nu^{\prime})\\
&=\sum_{l=1}^{n^{\prime}} \sum_{\substack{{x^{\prime}\in X},\\{\xi_1\in X^{i_1}},{\ldots},\\{\xi_l^{\prime}\in X^{i_l-1}},{\ldots},\\{\xi_{n^\prime}\in X^{i_{n^\prime}}}}}\delta_{{x}{x^{\prime}}}
(\xi_1\shuffle\cdots\shuffle\xi_l^{\prime}\shuffle\cdots\shuffle\xi_{n^\prime},\nu^{\prime})\\
&=\sum_{l=1}^{n^{\prime}} \sum_{\substack{{\xi_1\in X^{i_1}},{\ldots},\\{\xi_l^{\prime}\in X^{i_l-1}},{\ldots},\\{\xi_{n^\prime}\in X^{i_{n^\prime}}}}}
(\xi_1\shuffle\cdots\shuffle\xi_l^{\prime}\shuffle\cdots\shuffle\xi_{n^\prime},\nu^{\prime})\\
&=\sum_{l=1}^{n^{\prime}} \binom{j-1}{i_1, \ldots, (i_l-1), \ldots, i_{n^{\prime}}} = \binom{j}{i_1, \ldots, i_{n^{\prime}}},
\end{align*}
where $\delta_{xx^\prime}$ denotes the Kronecker delta symbol, and the multinomial coefficient recurrence formula was applied in the final step.
Hence, by induction, the equality holds for all $\nu\in X^*$.

\end{document}